\documentstyle[12pt]{amsart}
\begin{document}
\title{CR-rigidity of pseudo harmonic maps and pseudo biharmonic maps}
  \title[$CR$ rigidity of pseudo biharmonic maps]
  {$CR$ rigidity of pseudo harmonic maps and pseudo biharmonic maps}
   \author{Hajime Urakawa}
  \address{Division of Mathematics, Graduate School of Information Sciences, Tohoku University, Aoba 6-3-09, Sendai, 980-8579, Japan}
  \curraddr{Institute for International Education, 
  Tohoku University, Kawauchi 41, Sendai 980-8576, Japan}
  \email{urakawa@@math.is.tohoku.ac.jp}
  \thanks{
  Supported by the Grant-in-Aid for the Scientific Research, (C) No. 25400154, Japan Society for the Promotion of Science. 
  }
\maketitle
\begin{abstract} 
The $CR$ analogue of B.-Y. Chen's conjecture on pseudo biharmonic maps 
will be shown. Pseudo biharmonic, but not pseudo harmonic, 
isometric immersions with parallel pseudo mean curvature vector fields, will be characterized. 
    \end{abstract}
\numberwithin{equation}{section}
\theoremstyle{plain}
\newtheorem{df}{Definition}[section]
\newtheorem{th}[df]{Theorem}
\newtheorem{prop}[df]{Proposition}
\newtheorem{lem}[df]{Lemma}
\newtheorem{cor}[df]{Corollary}
\newtheorem{rem}[df]{Remark}
\section{Introduction}
Harmonic maps play a central role in geometry;\,they are critical points of the energy functional 
$E(\varphi)=\frac12\int_M\vert d\varphi\vert^2\,v_g$ 
for smooth maps $\varphi$ of $(M,g)$ into $(N,h)$. The Euler-Lagrange equations are given by the vanishing of the tension filed 
$\tau(\varphi)$. 
In 1983, Eells and Lemaire \cite{EL1} extended the notion of harmonic map to  
biharmonic map, which are 
critical points of the bienergy functional
$E_2(\varphi)=\frac12\int_M
\vert\tau(\varphi)\vert^2\,v_g$.
After Jiang \cite{J} studied the first and second variation formulas of $E_2$, 
extensive studies in this area have been done
(for instance, see 
\cite{CMP}, \cite{LO2},  \cite{MO1},
\cite{IIU2}, \cite{IIU},  \cite{II}). Every harmonic maps is always biharmonic by definition.  
Chen raised (\cite{C}) famous Chen's conjecture and later, Caddeo, Montaldo, Piu and Oniciuc raised (\cite{CMP}) the generalized Chen's conjecture. 
\vskip0.2cm\par
\textbf{B.-Y. Chen's conjecture:}\par
 {\em Every biharmonic submanifold of the Euclidean space ${\mathbb R}^n$ must be harmonic (minimal).}
\medskip\par
\textbf{The generalized B.-Y. Chen's conjecture:}\par
{\em Every biharmonic submanifold of a Riemannian manifold of non-positive curvature 
must be harmonic (minimal).}
\medskip\par
For the generalized Chen's conjecture, 
Ou and Tang gave (\cite{OT}) a counter example in a Riemannian manifold of negative curvature. 
For Chen's conjecture, some affirmative answers were known for 
surfaces in the three dimensional Euclidean space (\cite{C}), 
and 
hypersurfaces of the four dimensional Euclidean space (\cite{HV}, \cite{D}). 
Akutagawa and Maeta showed (\cite{AM}) that 
any complete regular biharmonic submanifold of the Euclidean space 
${\mathbb R}^n$ is harmonic (minimal). 
\par
To the generalized Chen's conjecture, we showed (\cite{NUG}) that: 
for a complete Riemannian manifold $(M,g)$, a Riemannian manifold 
$(N,h)$ of non-positive curvature, then,   
every biharmonic map 
$\varphi:\,(M,g)\rightarrow (N,h)$ with 
finite energy and 
finite bienergy 
is harmonic. 
In the case ${\rm Vol}(M,g)=\infty$, 
every biharmonic map 
$\varphi:\,(M,g)\rightarrow (N,h)$ with finite bienergy 
is harmonic. 
This 
gave (\cite{NU1}, \cite{NU2}, \cite{NUG}) affirmative answers to the generalized Chen's conjecture
under the $L^2$-condition and the completeness of $(M,g)$. 
\vskip0.2cm\par
In 1970's, Chern and Moser 
initiated (\cite{CM}) the geometry and analysis of strictly convex $CR$ manifolds, and many mathematicians works on $CR$ manifolds (cf. \cite{DT}). Recently, 
Barletta, Dragomir and Urakawa gave (\cite{BDU}) 
the notion of pseudo harmonic map, and also 
Dragomir and Montaldo settled (\cite{DM}) the one of pseudo biharmonic map. 
\par
In this paper, we raise
\vskip0.2cm\par 
{\bf The $CR$ analogue of the generalized Chen's conjecture:}
\par
{\em
Let $(M,g_{\theta})$ be a complete 
strictly pseudoconvex 
$CR$ manifold, and assume that $(N,h)$ is a Riemannian manifold of non-positive curvature. 
\par\noindent
Then, every pseudo biharmonic isometric immersion 
$\varphi:\,(M,g_{\theta})\rightarrow (N,h)$ must be 
pseudo harmonic.
}
\vskip0.2cm\par
We will show this conjecture holds under some $L^2$ condition on a complete strongly pseudoconvex $CR$ manifold (cf. Theorem 3.2), and 
will give characterization theorems on pseudo biharmonic immersions from 
$CR$ manifolds into the unit sphere or the complex projective space (cf. Theorems 6.2 and 7.1).  More precisely, we will show 
\begin{th}  (cf. Theorem 3.2) 
Let $\varphi$ be a pseudo biharmonic map of a complete $CR$ manifold $(M,g_{\theta})$ into a Riemannian manifold $(N,h)$ of non-positive curvature. Then, 
\par
If the pseudo energy $E_b(\varphi)$ and the pseudo bienergy $E_{b,2}(\varphi)$ are finite, then $\varphi$ is pseudo harmonic. 
\end{th}
\par
For isometric immersions of a $CR$ manifold $(M^{2n+1},g_{\theta})$ into the unit sphere $S^{2n+2}(1)$ of curvature $1$, 
we have 
\begin{th}  (cf. Theorem 6.2) 
For such immersion, assume that 
the pseudo mean curvature is parallel, but not pseudo harmonic. 
\par
Then, $\varphi$ is pseudo biharmonic if and only if 
the restriction of 
the second fundamental form $B_{\varphi}$ to 
the holomorphic subspace $H_x(M)$ of $T_xM$ $(x\in M)$ satisfies that 
$$\Vert\, B_{\varphi}\vert_{H(M)\times H(M)}\,\Vert^2=2n.$$ 
\end{th}
\par
For isometric immersions of a $CR$ manifold $(M^{2n+1},g_{\theta})$ into the complex projective space $({\mathbb P}^{n+1}(c),h,J)$ of holomorphic sectional curvature $c>0$, we have 
\begin{th}  (cf. Theorem 7.1) 
For such immersion, assume that 
the pseudo mean curvature is parallel, but not pseudo harmonic. 
Then, $\varphi$ is pseudo biharmonic if and only if 
one of the following holds: 
\par
$(1)$ $J(d\varphi(T))$ is tangent to $\varphi(M)$ and 
$$\Vert\, B_{\varphi}\vert_{H(M)\times H(M)}\,\Vert^2
=\frac{c}{4}(2n+3).$$
\par
$(2)$ $J(d\varphi(T))$ is normal to $\varphi(M)$ and 
$$\Vert\, B_{\varphi}\vert_{H(M)\times H(M)}\,\Vert^2
=\frac{c}{4}(2n)=\frac{n}{2}\,c.$$
Here, $T$ is the charactersitic vector field of $(M,g_{\theta})$, $H_x(M)\oplus {\mathbb R}T_x=T_x(M)$, and 
$B_{\varphi}\vert_{H(M)\times H(M)}$ is 
the restriction of 
the second fundamental form $B_{\varphi}$ to 
$H_x(M)$ $(x\in M)$. 
\end{th}
\par
Several examples of pseudo biharmonic immersions 
of $(M,g_{\theta})$ into the unit sphere or complex projective space will be given. 
\vskip0.2cm\par
{\bf Acknowledgement.} \quad 
This work was finished during the stay at the University of Basilicata, Potenza, Italy, September of 2014. 
The author was invited by Professor Sorin Dragomir to the University of Basilicata, Italy. 
The author 
would like to express his sincere gratitude to Professor Sorin Dragomir and 
Professor Elisabetta Barletta for their kind hospitality and helpful discussions.
\vskip0.6cm\par
\section{Preliminaries}
{\bf 2.1.}\quad We prepare the materials for the first and second variational formulas for the bienergy functional and biharmonic maps. 
Let us recall the definition of a harmonic map $\varphi:\,(M,g)\rightarrow (N,h)$, of a compact Riemannian manifold $(M,g)$ into another Riemannian manifold $(N,h)$, 
which is an extremal 
of the {\em energy functional} defined by 
$$
E(\varphi)=\int_Me(\varphi)\,v_g, 
$$
where $e(\varphi):=\frac12\vert d\varphi\vert^2$ is called the energy density 
of $\varphi$.  
That is, for any variation $\{\varphi_t\}$ of $\varphi$ with 
$\varphi_0=\varphi$, 
\begin{equation}
\frac{d}{dt}\bigg\vert_{t=0}E(\varphi_t)=-\int_Mh(\tau(\varphi),V)v_g=0,
\end{equation}
where $V\in \Gamma(\varphi^{-1}TN)$ is a variation vector field along $\varphi$ which is given by 
$V(x)=\frac{d}{dt}\big\vert_{t=0}\varphi_t(x)\in T_{\varphi(x)}N$, 
$(x\in M)$, 
and  the {\em tension field} is given by 
$\tau(\varphi)
=\sum_{i=1}^mB_{\varphi}(e_i,e_i)\in \Gamma(\varphi^{-1}TN)$, 
where 
$\{e_i\}_{i=1}^m$ is a locally defined orthonormal frame field on $(M,g)$, 
and $B_{\varphi}$ is the second fundamental form of $\varphi$ 
defined by 
\begin{align}
B_{\varphi}(X,Y)&=(\widetilde{\nabla}d\varphi)(X,Y)\nonumber\\
&=(\widetilde{\nabla}_Xd\varphi)(Y)\nonumber\\
&=\overline{\nabla}_X(d\varphi(Y))-d\varphi(\nabla^g_XY),
\end{align}
for all vector fields $X, Y\in {\frak X}(M)$. 
Here, 
$\nabla^g$, and
$\nabla^h$, 
 are Levi-Civita connections on $TM$, $TN$  of $(M,g)$, $(N,h)$, respectively, and 
$\overline{\nabla}$, and $\widetilde{\nabla}$ are the induced ones on $\varphi^{-1}TN$, and $T^{\ast}M\otimes \varphi^{-1}TN$, respectively. By (2.1), $\varphi$ is {\em harmonic} if and only if $\tau(\varphi)=0$. 
\par
The second variation formula is given as follows. Assume that 
$\varphi$ is harmonic. 
Then, 
\begin{equation}
\frac{d^2}{dt^2}\bigg\vert_{t=0}E(\varphi_t)
=\int_Mh(J(V),V)v_g, 
\end{equation}
where 
$J$ is an elliptic differential operator, called the
{\em Jacobi operator}  acting on 
$\Gamma(\varphi^{-1}TN)$ given by 
\begin{equation}
J(V)=\overline{\Delta}V-{\mathcal R}(V),
\end{equation}
where 
$\overline{\Delta}V=\overline{\nabla}^{\ast}\overline{\nabla}V
=-\sum_{i=1}^m\{
\overline{\nabla}_{e_i}\overline{\nabla}_{e_i}V-\overline{\nabla}_{\nabla^g_{e_i}e_i}V
\}$ 
is the {\em rough Laplacian} and 
${\mathcal R}$ is a linear operator on $\Gamma(\varphi^{-1}TN)$
given by 
${\mathcal R}(V)=
\sum_{i=1}^mR^h(V,d\varphi(e_i))d\varphi(e_i)$,
and $R^h$ is the curvature tensor of $(N,h)$ given by 
$R^h(U,V)=\nabla^h{}_U\nabla^h{}_V-\nabla^h{}_V\nabla^h{}_U-\nabla^h{}_{[U,V]}$ for $U,\,V\in {\frak X}(N)$.   
\par
J. Eells and L. Lemaire \cite{EL1} proposed polyharmonic ($k$-harmonic) maps and 
Jiang \cite{J} studied the first and second variation formulas of biharmonic maps. Let us consider the {\em bienergy functional} 
defined by 
\begin{equation}
E_2(\varphi)=\frac12\int_M\vert\tau(\varphi)\vert ^2v_g, 
\end{equation}
where 
$\vert V\vert^2=h(V,V)$, $V\in \Gamma(\varphi^{-1}TN)$.  
\par
The first variation formula of the bienergy functional 
is given by
\begin{equation}
\frac{d}{dt}\bigg\vert_{t=0}E_2(\varphi_t)
=-\int_Mh(\tau_2(\varphi),V)v_g.
\end{equation}
Here, 
\begin{equation}
\tau_2(\varphi)
:=J(\tau(\varphi))=\overline{\Delta}(\tau(\varphi))-{\mathcal R}(\tau(\varphi)),
\end{equation}
which is called the {\em bitension field} of $\varphi$, and 
$J$ is given in $(2.4)$.  
\par
A smooth map $\varphi$ of $(M,g)$ into $(N,h)$ is said to be 
{\em biharmonic} if 
$\tau_2(\varphi)=0$. 
By definition, every harmonic map is biharmonic. 
For an isometric immersion, it is minimal if and only if it is harmonic.  
\vskip0.6cm\par
{\bf 2.2.} \quad
Following Dragomir and Montaldo \cite{DM}, and also 
Barletta, Dragomir and Urakawa \cite{BDU}, we will prepare 
the materials on pseudo harmonic maps and pseudo biharmonic maps. 
\par
Let $M$ be a strictly pseudoconvex $CR$ manifold 
of $(2n+1)$-dimension, 
$T$, the characteristic vector field on $M$,  
$J$ is the complex structure of the subspace 
$H_x(M)$ of $T_x(M)$ $(x\in M)$, and 
$g_{\theta}$, the Webster Riemannian metric on $M$ defined for $X,Y\in H(M)$ by 
$$
g_{\theta}(X,Y)=(d\theta)(X,JY),\,\,
g_{\theta}(X,T)=0,\,\,g_{\theta}(T,T)=1.
$$  
\par
Let us recall for a $C^{\infty}$ map $\varphi$ of 
$(M,g_{\theta})$ into another Riemannian manifold $(N,h)$, 
the {\em pseudo energy} $E_b(\varphi)$ is defined (\cite{BDU}) by  
 \begin{equation}
 E_b(\varphi)=\frac12\int_M
 \sum_{i=1}^{2n}(\varphi^{\ast}h)
 (X_i,X_i)\,\theta\wedge (d\theta)^n,
 \end{equation}
 where $\{X_i\}_{i=1}^{2n}$ is an orthonormal frame field 
 on $(H(M),g_{\theta})$. 
 Then, the first variational formula of $E_b(\varphi)$ is as follows (\cite{BDU}). 
 For every variation $\{\varphi_t\}$ of $\varphi$ with 
 $\varphi_0=\varphi$, 
 \begin{align}
 \frac{d}{dt}\bigg\vert_{t=0}
 E_b(\varphi_t)
 =-\int_Mh(\tau_b(\varphi),V)\,d\theta\wedge(d\theta)^n=0,
 \end{align}
 where 
 $V\in \Gamma(\varphi^{-1}TN)$ is defined by 
 $V(x)=\frac{d}{dt}\vert_{t=0}\varphi_t(x)\in T_{\varphi(x)}N$, 
 $(x\in M)$. Here, $\tau_b(\varphi)$ is the {\em pseudo tension field} which is given by 
 \begin{equation}
 \tau_b(\varphi)=\sum_{i=1}^{2n}B_{\varphi}(X_i,X_i),
 \end{equation}
 where $B_{\varphi}(X,Y)$ $(X,\,Y\in {\frak X}(M)$) 
 is the second fundamental form (2.2) 
 for a $C^{\infty}$ map of $(M,g_{\theta})$ into $(N,h)$. 
 Then, $\varphi$ is {\em pseudo harmonic} if 
 $\tau_b(\varphi)=0$. 
 \par
 The second variational formula of $E_{b}$ is given as follows (\cite{BDU}, p.733): 
 \begin{align}
 \frac{d^2}{dt^2}\bigg\vert_{t=0}E_b(\varphi_t)
 =\int_Mh(J_b(V),V)\,\theta\wedge(d\theta)^n,
 \end{align}
 where 
 $J_b$ is a subelliptic operator acting on 
 $\Gamma(\varphi^{-1}TN)$ given by 
 \begin{align}
 J_b(V)=\Delta_b\,V-{\mathcal R}_b(V). 
 \end{align}
Here, for $V\in \Gamma(\varphi^{-1}TN))$,
\begin{equation}
\left\{
 \begin{aligned}
 \Delta_bV&=(\overline{\nabla}^H)^{\ast}\,\overline{\nabla}^H\,V
 =-\sum_{i=1}^{2n}\left\{
 \overline{\nabla}_{X_i}(\overline{\nabla}_{X_i}V)-
 \overline{\nabla}_{\nabla_{X_i}X_i}V
 \right\},\\
 {\mathcal R}_b(V)&=\sum_{i=1}^{2n}R^h(V,d\varphi(X_i))
 d\varphi(X_i),
 \end{aligned}
 \right.
 \end{equation}
 where $\nabla$ is the Tanaka-Webster connection, and
 $\overline{\nabla}$, the induced connection on 
 $\phi^{-1}TN$ induced from the Levi-Civita connection 
 $\nabla^h$, 
 and 
 $\{X_i\}_{i=1}^{2n}$, a local orthonormal frame field on 
 $(H(M),g_{\theta})$, respectively. 
 Here, $(\overline{\nabla}^H)_XV:=\overline{\nabla}_{X^H}V$ $(X\in {\frak X}(M), \,V\in \Gamma(\phi^{-1}TN))$, 
 corresponding to the decomposition 
 $X=X^H+g_{\theta}(X,T)\,T$ $(X^H\in H(M))$, 
 and define $\pi_H(X)=X^H$ $(X\in T_x(M))$, and 
 $(\overline{\nabla}^H)^{\ast}$ is 
 the formal adjoint of $\overline{\nabla}^H$.   
 \par
 Dragomir and Montaldo \cite{DM} introduced 
 the {\em pseudo bienergy} 
 given by 
 \begin{equation}
 E_{b,2}(\varphi)=\frac12\int_Mh(\tau_b(\varphi),\tau_b(\varphi))\,\theta\wedge (d\theta)^n, 
 \end{equation}
 where $\tau_b(\varphi)$ is the {\em pseudo tension field} 
 of $\varphi$. They gave the first variational formula 
 of $E_{b,2}$ as follows (\cite{DM}, p.227): 
 \begin{align}
 \frac{d}{dt}\bigg\vert_{t=0}E_{b,2}(\varphi_t)
= -\int_Mh(\tau_{b,2}(\varphi),V)\,\theta\wedge(d\theta)^n,
 \end{align} 
 where 
 $\tau_{b,2}(\varphi)$ is called the {\em pseudo bitension field} given by 
 \begin{align}
 \tau_{b,2}(\varphi)=\Delta_b\big(\tau_b(\varphi)\big)
 -\sum_{i=1}^{2n}R^h(\tau_b(\varphi),d\varphi(X_i))\,d\varphi(X_i). 
 \end{align}
 \par
 Then, a smooth map $\varphi$ of $(M,g_{\theta})$ into $(N,h)$ 
 is said to be {\em pseudo biharmonic} if 
 $\tau_{b,2}(\varphi)=0$. 
 By definition, 
 a pseudo harmonic map is always pseudo biharmonic. 
\medskip\par
\section{Generalized Chen's conjecture 
for pseudo biharmonic maps}
{\bf 3.1}\quad 
First, let us recall the usual  
Weitzenbeck formula 
for a $C^{\infty}$ map from a Riemannian manifod 
$(M,g)$ of $(2n+1)$ dimension 
 into a Riemannian manifold $(N,h)$: 
 \medskip\par
\begin{lem} (The Weitzenbeck formula)
\quad For every $C^{\infty}$ map 
$\varphi$ of $(M,g)$ 
of $(2n+1)$-dimension into a  Riemannian manifold $(N,h)$, 
the Hodge Laplacian $\Delta$ 
acting on the
$1$-form $d\varphi$, regarded as 
a $\varphi^{-1}TN$-valued $1$ form, 
$d\varphi\in \Gamma(T^{\ast}M\otimes \varphi^{-1}TN)$, we have 
\begin{align}
\Delta \,d\varphi=\widetilde{\nabla}^{\ast}\,\widetilde{\nabla}\,d\varphi+S.
\end{align}
Here, let us recall 
the rough Laplacian 
\begin{align}
\widetilde{\nabla}^{\ast}\,\widetilde{\nabla}
:=\sum_{k=1}^{2n+1}\left\{
\widetilde{\nabla}_{e_k}\,\widetilde{\nabla}_{e_k}-\widetilde{\nabla}_{\nabla^g{}_{e_k}e_k}
\right\}
\end{align}
\begin{align}
S(X):=-(\widetilde{R}(X,e_k)d\varphi)(e_k), 
\qquad (X\in {\frak X}(M)). 
\end{align}
Here,  
$\nabla^g$, $\nabla^h$  are the Levi-Civita connections of $(M,g)$, $(N,h)$, 
and 
$\widetilde{\nabla}$ is 
the induced connection 
on $T^{\ast}M\otimes \varphi^{-1}TN$ defined by 
$(\widetilde{\nabla}_Xd\varphi)(Y)=\overline{\nabla}_Xd\varphi(Y)-d\varphi(\nabla^g{}_XY)$, 
$\overline{\nabla}$ is the induced connection on 
$\varphi^{-1}TN$ given by 
$\overline{\nabla}_Xd\varphi(Y)=\nabla^h{}_{d\varphi(X)}d\varphi(Y)$, $(X,\,Y\in {\frak X}(M))$, 
and $\{e_k\}_{k=1}^{2n+1}$ 
is a locally defined orthonormal vector field 
on $(M,g)$.   
The curvature tensor field $\widetilde{R}$ in (3.3) 
is defined by 
\begin{align}
(\widetilde{R}(X,Y)d\varphi)(Z)
&:=
\overline{R}(X,Y)\,\,d\varphi(Z)-d\varphi(R^g(X,Y)Z)\nonumber\\
&\,\,=
R^h(d\varphi(X),d\varphi(Y))d\varphi(Z)-d\varphi(R^g(X,Y)Z), \nonumber
\end{align}
for $X,\,Y,\,Z\in {\frak X}(M)$, 
where 
$\overline{R}$, 
$R^g$, and $R^h$ are the curvature tensors of 
the induced connection $\overline{\nabla}$, 
$\nabla^g$ and $\nabla^h$, respectively. 
\end{lem} 
\medskip\par
Notice that 
for an isometric immersion 
$\varphi:\,(M,g)\rightarrow (N,h)$, it holds that 
\begin{equation}
(\widetilde{\nabla}_Xd\varphi)(Y)=B_{\varphi}(X,Y),
\qquad (X,\,Y\in {\frak X}(M)).
\end{equation}
{\bf 3.2} In this part, 
we first raise the $CR$ analogue of the generalized 
Chen's conjecture, and 
settle it 
for pseudo biharmonic maps with 
finite pseudo energy and finite pseudo bienergy. 
\par
Let us recall 
a strictly pseudoconvex $CR$ manifold (possibly non compact) $(M,g_{\theta})$ 
of $(2n+1)$-dimension, 
and the Webster Riemannian metric  $g_{\theta}$ given by 
$$
g_{\theta}(X,Y)=(d\theta)(X,JY),\quad 
g_{\theta}(X,T)=0,\quad 
g_{\theta}(T,T)=1
$$
for $X,\,Y\in H(M)$. 
Recall 
the material on the Levi-Civita 
connection $\nabla^{g_{\theta}}$ 
of $(M,g_{\theta})$. 
Due to Lemma 1.3, Page 38  in \cite{DT}, 
it holds that, 
\begin{align}
\nabla^{g_{\theta}}=\nabla+(\Omega-A)\otimes T+\tau\otimes\theta+2\,\theta\odot J, 
\end{align}
where $\nabla$ is the Tanaka-Webster connection, 
$\Omega=d\theta$, and 
$A(X,Y)=g_{\theta}(\tau X,Y)$, 
$\tau X=T_{\nabla}(T,X)$, and  $T_{\nabla}$ is the torsion tensor of $\nabla$. And also,  
$(\tau\otimes \theta)(X,Y)=\theta(Y)\,\tau X$, 
$(\theta\odot J)(X,Y)=\frac12\,\{\theta(X)\,JY+\theta(Y)\,JX\}$ 
for all vector fields $X,\,Y$ on $M$. 
Here, $J$ is the complex structure 
on $H(M)$ and is extended as an endomorphism 
on $(M)$ by $JT=0$.  
\par
Then, we have 
\begin{align}
\nabla^{g_{\theta}}_{X_k}X_k&=\nabla_{X_k}X_k-A(X_k,X_k)\,T,\\
\nabla^{g_{\theta}}_TT&=0,
\end{align}
where $\{X_k\}_{k=1}^{2n}$ is a locally defined  orthonormal frame field on $H(M)$ with respect to $g_{\theta}$, and $T$ is the characteristic vector field of $(M,g_{\theta})$. 
For (3.6), it follows from that $\Omega(X_k,X_k)=0$, $(\tau\otimes\theta)(X_k,X_k)=0$, and 
$(\theta\odot J)(X_k,X_k)=0$ since $\theta(X_k)=0$. 
For (3.7), notice that the Tanaka-Webster connection 
$\nabla$ satisfies $\nabla_TT=0$, and also $\tau T=0$ and $JT=0$, so that 
$\Omega(T,T)=0,\,A(T,T)=0, \,(\tau\otimes \theta)(T,T)=0\,(\theta\odot J)(T,T)=0$ which imply (3.7).   
\medskip\par
For (3.2) in the Weitenbeck formula in Lemma 3.1, 
 by taking 
 $\{X_k\,\,(k=1,\cdots,2n),\,\,T\}$,  
 as an orthonormal basis 
 $\{e_k\}$ of our $(M,g^{\theta})$, and due to (3.6) and (3.7), we have 
\begin{align}
(\widetilde{\nabla}^{\ast}\,\widetilde{\nabla}\,d\varphi)(X)&
=(\widetilde{\Delta}_b\,d\varphi)(X)\nonumber\\
&=-
\sum_{k=1}^{2n+1}\{
\widetilde{\nabla}_{e_k}\,\widetilde{\nabla}_{e_k}
-\widetilde{\nabla}_{\nabla^{g_{\theta}}_{e_k}e_k}\}
\,d\varphi(X)\nonumber\\
&=-\sum_{k=1}^{2n}\{
\widetilde{\nabla}_{X_k}\,\widetilde{\nabla}_{X_k}-\widetilde{\nabla}_{\nabla^{g_{\theta}}_{X_k}X_k}\}
\,d\varphi(X)
\nonumber\\
&\quad
-\{\widetilde{\nabla}_T\,\widetilde{\nabla}_T-\widetilde{\nabla}_{\nabla^{g_{\theta}}_TT}\}\,d\varphi(X)\nonumber\\
&=-\sum_{k=1}^{2n}\{
\widetilde{\nabla}_{X_k}\,\widetilde{\nabla}_{X_k}-\widetilde{\nabla}_{\nabla_{X_k}X_k}\}\,d\varphi(X)\nonumber\\
&\qquad-\{
\widetilde{\nabla}_T\,\widetilde{\nabla}_T+\sum_{k=1}^{2n}A(X_k,A_k)\,\widetilde{\nabla}_T\}\,d\varphi(X)\nonumber\\
&=-\sum_{k=1}^{2n}
\{
\widetilde{\nabla}_{X_k}\,\widetilde{\nabla}_{X_k}-\widetilde{\nabla}_{\nabla_{X_k}X_k}\}\,d\varphi(X)-
\widetilde{\nabla}_T\,\widetilde{\nabla}_T\,d\varphi(X).
\end{align}
since
$\sum_{k=1}^{2n}A(X_k,X_k)=0$ (cf. \cite{DT}, p. 35). 
\par
For (3.3) in the Weitzenbeck formula in Lemma 3.1, we have 
\begin{align}
S(X)&=-\sum_{k=1}^{2n+1}(\widetilde{R}(X,e_k)d\varphi)(e_k)\nonumber\\
&=-\sum_{k=1}^{2n}(\widetilde{R}(X,X_k)d\varphi)(X_k)-(\widetilde{R}(X,T)d\varphi)(T)\nonumber\\
&=-\sum_{k=1}^{2n}\left\{
R^h(d\varphi(X),d\varphi(X_k))d\varphi(X_k)-d\varphi(R^{g_{\theta}}(X,X_k)X_k)\right\}\nonumber\\
&\quad 
-\left\{
R^h(d\varphi(X),d\varphi(T))d\varphi(T)-d\varphi(R^{g_{\theta}}(X,T)T)\right\}.
\end{align}
And, we have the following formulas for (3.1) in our case,  
\begin{align}
\Delta\,d\varphi(X)&=d\,d^{\ast}\,d\varphi(X)
\nonumber\\
&=
-d\,\tau(\varphi)(X)
\nonumber\\
&=-\overline{\nabla}_X\tau(\varphi).
\end{align}
Therefore, we have 
\begin{align}
-(\widetilde{\Delta}_b\,d\varphi)(X)&=\sum_{k=1}^{2n}
\left\{
\widetilde{\nabla}_{X_k}\,\widetilde{\nabla}_{X_k}-\widetilde{\nabla}_{\nabla_{X_k}X_k}
\right\}\,d\varphi(X)\nonumber\\
&=-(\Delta\,d\varphi)(X)+S(X)-
\widetilde{\nabla}_T\,\widetilde{\nabla}_T
\,d\varphi(X)\nonumber\\
&=
\overline{\nabla}_X\tau(\varphi)\nonumber\\
&\quad 
-\sum_{k=1}^{2n}\{
R^h(d\varphi(X),d\varphi(X_k))d\varphi(X_k)
-d\varphi(R^{g_{\theta}}(X,X_k)X_k)\}\nonumber\\
&\quad
-\{R^h(d\varphi(X),d\varphi(T))d\varphi(T)
-d\varphi(R^{g_{\theta}}(X,T)T)\}
\nonumber\\
&\quad 
-\widetilde{\nabla}_T\widetilde{\nabla}_Td\varphi(X).
\end{align}
\medskip\par
{\bf 3.3}\quad 
Let us consider 
the generalized B.-Y. Chen's conjecture 
for pseudo biharmonic maps which is $CR$ analogue of the usual generalized Chen's conjecture for biharmonic maps:
\vskip0.2cm\par
{\bf The $CR$ analogue of the generalized B.-Y. Chen's conjecture 
for pseudo biharmonic maps:}
\par
{\em
Let $(M,g_{\theta})$ be a complete 
strictly pseudoconvex 
$CR$ manifold, and assume that $(N,h)$ is a Riemannian manifold of non-positive curvature. 
\par
Then, every pseudo biharmonic isometric immersion 
$\varphi:\,(M,g_{\theta})\rightarrow (N,h)$ must be 
pseudo harmonic.
}
\vskip0.6cm\par
In this section, we want to show that the above conjecture is true under the finiteness of the 
pseudo energy and pseudo bienergy. 
\vskip0.2cm\par
\begin{th}
Assume that $\varphi$ is a pseudo biharmonic map  
of a strictly pseudoconvex complete $CR$ manifold $(M,g_{\theta})$ into another Riemannian manifold $(N,h)$ of non positive curvature.  
\par
If $\varphi$ has finite pseudo bienergy $E_{b,2}(\varphi)<\infty$ and 
finite pseudo energy $E_b(\varphi)<\infty$, then it is  
pseudo harmonic, i.e., 
$\tau_b(\varphi)=0$. 
\end{th}
\medskip\par
(\textit{Proof of Theorem 3.2})\quad 
The proof is divided into several steps. 
\medskip\par
(The first step)\quad 
For an arbitrarily fixed point $x_0\in M$, let 
$B_r(x_0)=\{x\in M:\,r(x)<r\}$ where 
$r(x)$ is a distance function on $(M,g_{\theta})$, and 
let us take a cut off function $\eta$ on $(M,g_{\theta})$, i.e., 
\begin{equation}
\left\{
\begin{aligned}
&0 \leq \eta(x)\leq 1
\qquad\qquad (x\in M),\\
&\eta(x)=1
\qquad\qquad (x\in B_r(x_0)),\\
&\eta(x)=0
\qquad\qquad (x\not\in B_{2r}(x_0)),\\
&\vert\nabla^{g_{\theta}}\,\eta\vert
\leq \frac{2}{r}
\qquad\qquad (x\in M),
\end{aligned}
\right.
\end{equation} 
where $r$, $\nabla^{g_{\theta}}$ 
are the distance function, 
the Levi-Civita connection of 
$(M,g_{\theta})$, respectively. 
Assume that 
$\varphi:\,(M,g_{\theta})\rightarrow (N,h)$ is a {\em pseudo biharmonic map}, i.e., 
\begin{align}
\tau_{b,\,2}(\varphi)&=J_b(\tau_b(\varphi))\nonumber\\
&=\Delta_b(\tau_b(\varphi))-\sum_{j=1}^{2n}R^h(\tau_b(\varphi),d\varphi(X_j))\,d\varphi(X_j)\nonumber\\
&=0.
\end{align}
\par
(The second step)\quad 
Then, we have 
\begin{align}
\int_M&\langle \Delta_b(\tau_b(\varphi)),\eta^2\,\tau_b(\varphi)\rangle\,\theta\wedge(d\theta)^n\nonumber\\
&=\int_M
\eta^2\,\sum_{j=1}^{2n}\langle R^h(\tau_b(\varphi),d\varphi(X_j))\,d\varphi(X_j),\,\tau_b(\varphi)\rangle\,\theta\wedge 
(d\theta)^n\nonumber\\
&\leq 0
\end{align}
since 
$(N,h)$ has the non-positive sectional curvature. 
But, for the left hand side of (3.14), it holds that 
\begin{align}
\int_M&\langle \Delta_b(\tau_b(\varphi)),\eta^2\,\tau_b(\varphi)\rangle\,\theta\wedge(d\theta)^n\nonumber\\
&=\int_M
\langle
\overline{\nabla}^H\,\tau_b(\varphi),\overline{\nabla}^H\,(\,\eta^2\,\tau_b(\varphi)\,)\,\rangle\,\theta\wedge(d\theta)^n\nonumber\\
&=\int_M
\sum_{j=1}^{2n}\langle
\overline{\nabla}_{X_j}\tau_b(\varphi),\overline{\nabla}_{X_j}(\eta^2\,\tau_b(\varphi))\,\rangle\,\theta\wedge(d\theta)^n.
\end{align}
Here, let us recall, for $V,\,W\in \Gamma(\varphi^{-1}TN))$, 
$$\langle\overline{\nabla}^HV,\overline{\nabla}^HW\rangle=\sum_{\alpha}\langle\overline{\nabla}^H_{e_{\alpha}}V,\overline{\nabla}^H_{e_{\alpha}}W\rangle
=\sum_{j=1}^{2n}\langle
\overline{\nabla}_{X_i}V,\overline{\nabla}_{X_i}W\rangle,$$
where 
$\{e_{\alpha}\}$ is a locally defined orthonormal frame field of $(M,g_{\theta})$ and 
$\overline{\nabla}^H_XW$ $(X\in {\frak X}(M),\,W\in \Gamma(\varphi^{-1}TN))$ is defined by 
$$\overline{\nabla}^H_XW=\sum_j\,\{
(X^Hf_j)\,V_j+f_j\overline{\nabla}_{X^H}\,V_j\}$$
for 
$W=\sum_jf_i\,V_j$ 
$(f_j\in C^{\infty}(M)$ and $V_j\in \Gamma(\varphi^{-1}TN)$. 
Here, 
 $X^H$ is the $H(M)$-component of $X$ 
 corresponding to the decomposition of 
 $T_x(M)=H_x(M)\oplus {\mathbb R}T_x$ $(x\in M)$, 
 and $\overline{\nabla}$ is the induced connection of
 $\varphi^{-1}TN$ from the Levi-Civita connection 
 $\nabla^h$ 
 of $(N,h)$. 
 \par
 Since 
 \begin{align}
 \overline{\nabla}_{X_j}(\eta^2\,\tau_b(\varphi))=
 2\eta\,X_j\eta\,\tau_b(\varphi)+\eta^2\,\overline{\nabla}_{X_j}\tau_b(\varphi), 
 \end{align}
 the right hand side of (3.15) is equal to 
 \begin{align}
 \int_M&\eta^2\,\sum_{j=1}^{2n}\big\vert\,\overline{\nabla}_{X_j}\tau_b(\varphi)\,\big\vert^2\,\theta\wedge(d\theta)^n\nonumber\\
 &+2\int_M\sum_{j=1}^{2n}\langle\eta\,\overline{\nabla}_{X_j}\tau_b(\varphi),\,(X_j\eta)\,\tau_b(\varphi)\rangle\,\theta\wedge(d\theta)^n. 
 \end{align}
 Therefore, together with (3.14), we have 
 \begin{align}
 \int_M&\eta^2\,\sum_{j=1}^{2n}\big\vert\,\overline{\nabla}_{X_j}\,\tau_b(\varphi)\,\big\vert^2\,\theta\wedge(d\theta)^n\nonumber\\
 &\leq -2\int_M\sum_{j=1}^{2n}\langle\eta\,\overline{\nabla}_{X_j}\tau_b(\varphi),\,(X_j\eta)\,\tau_b(\varphi)\rangle\,\theta\wedge(d\theta)^n\nonumber\\
 &=:-2\int_M\sum_{j=1}^{2n}\langle V_j,W_j\rangle\,\theta\wedge(d\theta)^n, 
 \end{align}
 where we define $V_j,\,W_j\in \Gamma(\varphi^{-1}TN)$ $(j=1,\cdots,2n)$ by 
 $$
 V_j:=\eta\,\overline{\nabla}_{X_j}\,\tau_b(\varphi),\quad 
 W_j:=(X_j\eta)\,\tau_b(\varphi).
 $$
 Then, since it holds that 
 $0\leq \big\vert \sqrt{\epsilon}\,V_i\pm\frac{1}{\sqrt{\epsilon}}\,W_i\big\vert^2$ for every $\epsilon>0$, we have,  
 \begin{align}
 &\mbox{\rm the right hand side of (3.18)}\nonumber\\
 &\qquad\leq 
 \epsilon\int_M\sum_{j=1}^{2n}\big\vert \,V_j\,\big\vert^2\,\theta\wedge(d\theta)^n+\frac{1}{\epsilon}\int_M\sum_{j=1}^{2n}\big\vert\,W_j\,\big\vert^2\,\theta\wedge(d\theta)^n
 \end{align}
 foe every $\epsilon>0$. By taking $\displaystyle\epsilon=\frac12$, we obtain 
 \begin{align}
 &\int_M\eta^2\sum_{j=1}^{2n}\big\vert\,\overline{\nabla}_{X_j}\tau_b(\varphi)\,\big\vert^2\,\theta\wedge(d\theta)^n\nonumber\\
 &\leq
 \frac12\int_M\sum_{j=1}^{2n}\eta^2\,\big\vert\overline{\nabla}_{X_j}\tau_b(\varphi)\big\vert^2\,\theta\wedge(d\theta)^n+
 2\int_M\sum_{j=1}^{2n}\big\vert X_j\eta\big\vert^2\big\vert\tau_b(\varphi)\,\big\vert^2\,\theta\wedge(d\theta)^n. 
 \end{align}
 Therefore, we obtain, due to the properties that 
 $\eta=1$ on $B_r(x_0)$, and 
 $\sum_{j=1}^{2n}\vert X_j\eta\vert^2\leq \vert\nabla^{g_{\theta}}\eta\vert^2\leq \left(\frac{2}{r}\right)^2$,
 \begin{align}
 \int_{B_r(x_0)}\sum_{j=1}^{2n}\big\vert\overline{\nabla}\tau_b(\varphi)\big\vert^2\,\theta\wedge(d\theta)^n
 &\leq 
 \int_M\eta^2\sum_{j=1}^{2n}\big\vert\overline{\nabla}_{X_j}\tau_b(\varphi)\big\vert^2\,\theta\wedge(d\theta)^n\nonumber\\
 &\leq 4\int_M\sum_{j=1}^{2n}\big\vert X_j\eta\big\vert^2\,\big\vert\tau_b(\varphi)\big\vert^2\,\theta\wedge(d\theta)^n\nonumber\\
 &\leq \frac{16}{r^2}\int_M\vert\tau_b(\varphi)\vert^2\,\theta\wedge(d\theta)^n.
 \end{align}
 \medskip\par
 (The third step)\quad By our assumption that 
 $E_{b,\,2}(\varphi)=\frac12\int_M\vert\tau_b(\varphi)\vert^2\,\theta\wedge(d\theta)^n<\infty$ and $(M,g_{\theta})$ is complete, 
 if we let $r\rightarrow \infty$, then $B_r(x_0)$ goes to $M$, and 
 the right hand side of (3.21) goes to zero. We have 
 \begin{align}
 \int_M\sum_{j=1}^{2n}\big\vert\,\overline{\nabla}_{X_j}\tau_b(\varphi)\,\big\vert^2\,\theta\wedge(d\theta)^n=0.
 \end{align} 
 This implies that 
 \begin{align}
 \overline{\nabla}_X\tau_b(\varphi)=0\qquad\qquad (\mbox{for all}\,\,X\in H(M)).
 \end{align}
 \medskip\par
(The forth step)\quad  Let us take 
a $1$ form $\alpha$ on $M$ defined by 
$$\alpha(X)=
\left\{
\begin{aligned}
\langle d\varphi(X),&\tau_b(\varphi)\rangle, 
\qquad 
(X\in H(M)),\\
&0\qquad\qquad\,\,\, (X=T).
 \end{aligned}
 \right.
 $$
   Then, we have 
\begin{align}
\int_M\vert\alpha\vert\,\theta\wedge(d\theta)^n&=\int_M\left(
\sum_{j=1}^{2n}\alpha(X_j)\vert^2
\right)^{\frac12}\,\theta\wedge(d\theta)^n\nonumber\\
&\leq 
\left(\vert d_b\varphi\vert^2\,\theta\wedge(d\theta)^n
\right)^{\frac12}\,\left(
\int_M\vert\tau_b(\varphi)\vert^2\,\theta\wedge(d\theta)^n\right)^{\frac12}\nonumber\\
&=2\,\sqrt{E_b(\varphi)\,E_{b,2}(\varphi)}<\infty, 
\end{align}
where we put 
$d_b\varphi:=\sum_{i=1}^{2n}d\varphi(X_i)\otimes X_i$, 
\begin{equation}
\vert d_b\varphi\vert^2
=\sum_{i,j=1}^{2n}g_{\theta}(X_i,X_j)\,
h(d\varphi(X_i),d\varphi(X_j))=\sum_{i=1}^{2n}h(d\varphi(X_i),d\varphi(X_i)), \nonumber
\end{equation}
and 
\begin{align}
E_b(\varphi)=\frac12\int_M\vert d_b\varphi\vert^2\,\,\theta\wedge(d\theta)^n.
\end{align}
Furthermore, let us define a $C^{\infty}$ function 
$\delta_b\alpha$ on $M$ by 
\begin{align}
\delta_b\alpha
=-\sum_{j=1}^{2n}({\nabla}_{X_j}\alpha)(X_j)
=-\sum_{j=1}^{2n}\left\{
X_j(\alpha(X_j))-\alpha(\nabla_{X_j}X_j)
\right\}, 
\end{align}
where $\nabla$ is the Tanaka-Webster connection. 
Notice that 
\begin{align}
\mbox{\rm div}(\alpha)&=\sum_{j=1}^{2n}(\nabla^{g_{\theta}}_{X_j}\alpha)(X_j)+(\nabla^{g_{\theta}}_T\alpha)(T)\nonumber\\
&=\sum_{j=1}^{2n}\big\{
X_j(\alpha\,\circ\,\pi_H(X_j))-\alpha\,\circ\,\pi_H(\nabla^{g_{\theta}}_{X_j}X_j)\big\}\nonumber\\
&\qquad +T(\alpha\,\circ\,\pi_H(T))-\alpha\,\circ\,\pi_H(\nabla^{g_{\theta}}_TT)\nonumber\\
&=\sum_{j=1}^{2n}\big\{
X_j(\alpha(X_j))-\alpha(\pi_H(\nabla^{g_{\theta}}_{X_j}X_j))\big\}\nonumber\\
&=\sum_{j=1}^{2n}\big\{
X_j(\alpha(X_j))-\alpha(\nabla_{X_j}X_j)
\big\}\nonumber\\
&=-\delta_b\alpha,
\end{align}
where $\pi_H:\,T_x(M)\rightarrow H_x(M)$ is the natural projection. 
We used the facts that 
$\nabla^{g_{\theta}}_TT=0$, and $\pi_H(\nabla^{g_{\theta}}_{X}Y)=\nabla_{X}Y$ $(X,Y\in H(M))$ (\cite{A}, p.37). 
Here, recall again 
$\nabla^{g_{\theta}}$ is the Levi-Civita connection of 
$g_{\theta}$, and $\nabla$ is the Tanaka-Webster connection. 
Then, we have, for (3.26),  
\begin{align}
\delta_b\alpha&=-\sum_{j=1}^{2n}
\left\{X_j\,\langle d\varphi(X_j),\tau_b(\varphi)\rangle-
\langle d\varphi(\nabla_{X_j}X_j),\tau_b(\varphi)\rangle
\right\}\nonumber\\
&=-\sum_{j=1}^{2n}\bigg\{
\begin{aligned}
&\langle\overline{\nabla}_{X_j}(d\varphi(X_j)),\tau_b(\varphi)\rangle+
\langle\,d\varphi(X_j),\overline{\nabla}_{X_j}\tau_b(\varphi)\rangle\nonumber\\
&-\langle\,d\varphi(\nabla_{X_j}X_j),\tau_b(\varphi)\rangle
\end{aligned}
\bigg\}\nonumber\\
&=-\left\langle
\sum_{j=1}^{2n}\left\{\overline{\nabla}_{X_j}(d\varphi(X_j))-d\varphi(\nabla_{X_j}X_j)
\right\},\tau_b(\varphi)\right\rangle\nonumber\\
&=-\vert\tau_b(\varphi)\vert^2.
\end{align}
We used (3.23) $\overline{\nabla}_{X_j}\tau_b(\varphi)=0$  
to derive the last second equality of (3.28). 
Then, due to (3.28), we have 
for $E_{b,2}(\varphi)$, 
\begin{align}
E_{b,2}(\varphi)&=\frac12\int_M\vert\tau_b(\varphi)\vert^2\,\theta\wedge(d\theta)^n\nonumber\\
&=-\frac12\int_M\delta_b\alpha\,\,\theta\wedge(d\theta)^n\nonumber\\
&=\frac12\int_M\mbox{\rm div}(\alpha)\,\theta\wedge(d\theta)^n\nonumber\\
&=0.
\end{align} 
In the last equality, we used Gaffney's theorem (\cite{NUG}, p. 271, \cite{G}). 
\par
Therefore, we obtain $\tau_b(\varphi)\equiv 0$, i.e., 
$\varphi$ is pseudo harmonic. 
\qed
\vskip0.6cm\par
 \section{Parallel pseudo biharmonic isometric immersion into rank one symmetric spaces}
 On the contrary of the Section Three, 
 we consider isometric immersions into the unit sphere 
 or the complex projective spaces 
 which are pseudo biharmonic.   
One of the main theorem of this section is as follows: 
\vskip0.2cm\par
\begin{th}
Let $\varphi:\,(M,g_{\theta})\rightarrow S^{2n+2}(1)$ 
be an isometric immersion of 
a $CR$ manifold $(M,g_{\theta})$ of $(2n+1)$-dimension into the unit sphere $S^{2n+2}(1)$ of constant sectional curvature $1$ and $(2n+2)$-dimension. 
Assume that 
$\varphi$ admits a parallel pseudo mean curvature vector field with non-zero pseudo mean curvature. 
The following equivalences hold: 
\quad The immersion 
$\varphi$ is pseudo biharmonic if and only if 
\begin{equation}
\sum_{i=1}^{2n}\lambda_i{}^2=2n
\end{equation} 
if and only if 
\begin{equation}
\big\Vert B_{\varphi}\big\vert_{H(M)\times H(M)}\big\Vert^2=2n,
\end{equation}
where 
$\lambda_i$ $(1\leq i\leq 2n+1$ are the principal curvatures of the immersion $\varphi$ whose 
$\lambda_{2n+1}$ corresponds to the characteristic vector field $T$ of $(M,g_{\theta})$, and 
$B_{\varphi}\vert_{H(M)\times H(M)}$ is the restriction of the second fundamental form 
od $\varphi$ to the orthogonal complement $H(M)$ of $T$ in the tangent space $(T_x(M),g_{\theta})$.   
\end{th}
\vskip0.6cm\par
As applications of this theorem, 
we will give 
pseudo biharmonic immersions 
into the unit sphere 
which are not pseudo harmonic. 
  \vskip0.6cm\par
  The other case of rank one symmetric space is the complex projective space 
  ${\mathbb P}^{n+1}(c)$. 
  We obtain the following theorem: 
  \vskip0.2cm\par
  \begin{th}
  Let $\varphi:\,(M^{2n+1},g_{\theta} \rightarrow 
  {\mathbb P}^{n+1}(c)$ be an isometric immersion of 
  $CR$ manifold $(M,g_{\theta})$ into 
  the complex projective space ${\mathbb P}^{n+1}(c)$ of constant holomorphic sectional curvature $c$ and  complex $(n+1)$-dimension. 
  Assume that $\varphi$ has parallel pseudo-mean curvature vector filed with non-zero pseudo mean curvature. Then, the following equivalence relation holds: 
 The immersion $\varphi$ is pseudo-biharmonic if and only if 
 the following hold: 
 \par
 Either 
 $(1)$ $J(d\varphi(T))$ is tangent to 
 $\varphi(M)$ and 
 \begin{equation}\big\Vert B_{\varphi}\big\vert_{H(M)\times H(M)}\big\Vert^2=\frac{c}{4} \,(2n+3),
 \end{equation} 
 or 
 $(2)$ $J(d\varphi(T))$ is normal to $\varphi(M)$ and 
 \begin{equation}
 \big\Vert B_{\varphi}\big\vert_{H(M)\times H(M)}\big\Vert^2=\frac{c}{4}\,(2n)=c\,\frac{n}{2}.
 \end{equation}
  \end{th}
  \vskip0.6cm\par
  As applications of this theorem, we will give 
 pseudo biharmonic, but not pseudo harmonic immersions 
 $(M,g_{\theta})$ into the complex projective space 
 ${\mathbb P}^{n+1}(c)$.    
\section{Admissible immersions of strongly pseudoconvex $CR$ manifolds}
In this section, we introduce the notion of admissible isometric immersion of strongly pseudoconvex $CR$ manifold $(M,g_{\theta})$, and will show the following two lemmas related to $\Delta_b(\tau_b(\varphi))$ 
which are necessary to prove main theorems. 
\begin{df}
Let $(M^{2n+1},g_{\theta})$ be a strictly pseudoconvex $CR$ manifold, and 
$T_xM=H_x(M)\oplus {\mathbb R}T_x$, $(x\in M)$, 
the orthogonal decomposition of the tangent space $T_xM$ $(x\in M)$, where $T$ is the characteristic vector field of $(M^{2n+1},g_{\theta})$,   
$\varphi:\,(M^{2n+1},g_{\theta})\rightarrow (N,h)$ be an 
isometric immersion. 
The immersion $\varphi$ is called to be {\em admissible} if 
the second fundamental form $B_{\varphi}$ satisfies that 
\begin{equation}
B_{\varphi}(X,T)=0
\end{equation}
for all vector field $X$ in $H(M)$. 
\end{df}
\medskip\par
The following clarifies the meaning 
of the admissibility condition: 
\begin{prop}
Let $\varphi$ be an isometric immersion of a strongly pseudoconvex $CR$ manifold $(M^{2n+1},g_{\theta})$ into another Riemannian manifold $(N,h)$. 
Then, $\varphi$ is admissible if and only if 
\par
$(1)$ 
$d\varphi(T_x)$ $(x\in M)$ is a principal curvature 
vector field along $\varphi$ with some principal curvature 
$\lambda(x)$ $(x\in M)$. 
\par
This is equivalent the following:   
\par
$(2)$ 
The shape operator 
$A_{\xi}$ of the immersion $\varphi:\,(M,g_{\theta})\rightarrow (N,h)$ 
preserves $H_x(M)$ $(x\in M)$ invariantly 
for a normal vector field $\xi$. 
\end{prop}
\medskip\par
({\em Proof of Proposition 5.2})\quad 
We first note for every normal vector field $\xi$ of the isometric immersion 
$\varphi:\,(M,g_{\theta})\rightarrow (N,h)$, it holds that 
\begin{align*}
\quad
\langle B_{\varphi}(X,T),\xi\rangle=g_{\theta}(A_{\xi}X,T)
=g_{\theta}(X,A_{\xi}T),\qquad (X\in H_x(M)).\quad(\#) 
\end{align*}
Thus, if $\varphi$ is admissible, then 
the left hand side of $(\#)$ vanishes, then we have immediately that 
\begin{equation*}\qquad \quad 
\left\{
\begin{aligned}
&A_{\xi}X\in H_x(M)\qquad (X\in H_x(M)),\\
&A_{\xi}T_x=\lambda(x)\,T_x\quad (\mbox{for some real number 
$\lambda(x)$}).
\end{aligned}
\right.\qquad \qquad\quad
(\flat)
\end{equation*}
\par
Conversely, if one of the conditions of $(\flat)$ holds, 
then it turns out immediately that 
$\varphi$ is admissible. 
\qed
\bigskip\par
The following two lemmas will be essential to us 
later. 
\begin{lem}
Let $\varphi:\,(M^{2n+1},g_{\theta})\rightarrow (N,h)$ be an 
admissible 
isometric immersion 
with parallel pseudo mean curvature vector field. 
Then, 
the pseudo tension field $\tau_b(\varphi)$ 
satisfies that 
\begin{align}
-\Delta_b(\tau_b(\varphi))&=
\left\langle -\Delta_b(\tau_b(\varphi)),d\varphi(X_i)\right\rangle\,d\varphi(X_i)\nonumber
\\
&\quad 
+\left\langle \overline{\nabla}_{X_i}\tau_b(\varphi),d\varphi(X_j)\right\rangle\,
\left(\widetilde{\nabla}_{X_i}d\varphi\right)(X_j),
\end{align}  
where 
$\{X_j\}_{j=1}^{2n}$ is a local orthonormal frame field 
of $H(M)$ with respect to $g_{\theta}$. 
\end{lem}
\medskip\par
\begin{lem}
Under the same assumptions of the above lemma, 
we have 
\begin{align}
-\Delta_b(\tau_b(\varphi))&=
\left\langle \tau_b(\varphi), 
R^h(d\varphi(X_j),d\varphi(X_k))d\varphi(X_k)
\right\rangle\,d\varphi(X_j)\nonumber\\
&\quad 
+\left\langle \tau_b(\varphi),R^h(d\varphi(X_j),d\varphi(T))d\varphi(T)\right\rangle\,d\varphi(X_j)\nonumber\\
&\quad -\left\langle\tau_b(\varphi),\left(\widetilde{\nabla}_{X_i}d\varphi\right)(X_j)\right\rangle\,
\left(\widetilde{\nabla}_{X_i}d\varphi\right)(X_j),
\end{align}
where 
$R^h(U,V)W$ is the curvature tensor field of $(N,h)$ defined by 
$\displaystyle R^h(U,V)W=\nabla^h_U(\nabla^h_VW)-
\nabla^h_V(\nabla^N_UW)-\nabla^h_{[U,V]}W$ 
for vector fields $U,\,V,\,W$ on $N$, and 
$\nabla^h$ is the Levi-Civita connection of $(N,h)$.  
\end{lem}
\vskip0.6cm\par
({\em Proof of Lemma 5.3})     \quad 
The proof is divided into several steps. 
\par
(The first step)\quad 
Since we assume 
the pseudo mean curvature vector field 
$\tau_b(\varphi)$ is parallel, i.e., 
$\overline{\nabla}_X^{\perp}\tau_b(\varphi)=0$ 
$(X\in {\frak X}(M))$, 
the induced connection $\overline{\nabla}$ of the Levi-Civita connection $\nabla^h$ to the induced bundle 
$\varphi^{-1}TN$ satisfies that,   
for all $X\in {\frak X}(M)$,   
$$\overline{\nabla}_X\tau_b(\varphi)=\overline{\nabla}^{\top}_X\tau_b(\varphi)+\overline{\nabla}^{\perp}_X\tau_b(\varphi)=\overline{\nabla}^{\top}_X\tau_b(\varphi)\in \Gamma(\varphi_{\ast}TM).
$$
Then. we have, for all $X\in H(M)$, 
\begin{align}
\overline{\nabla}_X \tau_b(\varphi)
&=
\sum_{j=1}^{2n}\langle \overline{\nabla}_X\tau_b(\varphi),d\varphi(X_j)\rangle\,d\varphi(X_j)
+\langle \overline{\nabla}_X\tau_b(\varphi),d\varphi(T)\rangle\,d\varphi(T)\nonumber\\
&=\sum_{j=1}^{2n}\langle \overline{\nabla}_X\tau_b(\varphi),d\varphi(X_j)\rangle\,d\varphi(X_j). 
\end{align}
Due to the assumption of the admissibility of $\varphi$, 
for all $X\in H(M)$, 
\begin{align}
\langle\overline{\nabla}_X\tau_b(\varphi),d\varphi(T)\rangle 
=X\,\langle\tau_b(\varphi),d\varphi(T)\rangle-
\langle \tau_b(\varphi),\overline{\nabla}_Xd\varphi(T)\rangle=0.
\end{align}
In fact, 
$\tau_b(\varphi)=\sum_{i=1}^{2n}B_{\varphi}(X_i,X_i)$ is orthogonal to $d\varphi(TM)$ with respect to $\langle\,\,,\,\,\rangle$, 
we have $\langle\tau_b(\varphi),d\varphi(T)\rangle=0$. So, 
the first term of  (5.5) vanishes. 
By the admissibility of $\varphi$, 
for all $X\in H(M)$, 
\begin{equation}
0=B_{\varphi}(X,T)=\overline{\nabla}_Xd\varphi(T)-d\varphi(\nabla^{g_{\theta}}_XT),
\end{equation}
so that $\overline{\nabla}_Xd\varphi(T)$ is tangential, which implies that 
$$
\langle\tau_b(\varphi),\overline{\nabla}_Xd\varphi(T)\rangle=0.
$$
 We have (5.5), and then (5.4). 
 \par
 (The second step)\quad We calculate 
 $-\Delta_b(\tau_b(\varphi))$. 
 We have by (5.4),
 \begin{align}
 -\Delta_b(\tau_b(\varphi))&=
 \sum_{i=1}^{2n}
 \left\{
 \overline{\nabla}_{X_i}(\overline{\nabla}_{X_i}\tau_b(\varphi))-\overline{\nabla}_{\nabla_{X_i}X_i}\tau_b(\varphi)
 \right\}\nonumber\\
 &=\sum_{i=1}^{2n}
 \left[
 \begin{aligned}
 &\sum_{j=1}^{2n}
 \overline{\nabla}_{X_i}\left\{
 \langle \overline{\nabla}_X\tau_b(\varphi),d\varphi(X_j)\rangle\,d\varphi(X_j)
 \right\}\\
 &-\sum_{j=1}^{2n}
 \langle
 \overline{\nabla}_{\nabla_{X_i}X_i}\tau_b(\varphi),
 d\varphi(X_j)\rangle\,
 d\varphi(X_j)
 \end{aligned}
 \right]
\nonumber \\
 &=\sum_{i,j=1}^{2n}
 \left[
 \begin{aligned}
 &\langle
 \overline{\nabla}_{X_i}(\overline{\nabla}_{X_i}\tau_b(\varphi)),d\varphi(X_j)\rangle\,d\varphi(X_j)\\
 &+\langle\overline{\nabla}_{X_i}\tau_b(\varphi),\overline{\nabla}_{X_i}(d\varphi(X_j))\rangle\,d\varphi(X_j)\\
 &+\langle\overline{\nabla}_{X_i}\tau_b(\varphi),d\varphi(X_j)\rangle\,\overline{\nabla}_{X_i}d\varphi(X_j)\\
 &-\langle
 \overline{\nabla}_{\nabla_{X_i}X_i}\tau_b(\varphi),d\varphi(X_j)\rangle\,d\varphi(X_j)
 \end{aligned}
 \right]
 \nonumber\\
 &=
 \sum_{j=1}^{2n}\langle-\Delta_b(\tau_b(\varphi)),d\varphi(X_j)\rangle\,d\varphi(X_j)\nonumber\\
 &\quad +\sum_{i,j=1}^{2n}
 \left[
 \begin{aligned}
& \langle\overline{\nabla}_{X_i}\tau_b(\varphi),\overline{\nabla}_{X_i}(d\varphi(X_j))\rangle\,d\varphi(X_j)\\
 &+\langle \overline{\nabla}_{X_i}\tau_b(\varphi),d\varphi(X_j)\rangle\,\overline{\nabla}_{X_i}d\varphi(X_j)
 \end{aligned}
 \right].
 \end{align}
 (The third step)\quad 
 Here, 
 we have 
 \begin{equation}
 \left\{
 \begin{aligned}
& (\widetilde{\nabla}_{X_i}d\varphi)(X_j)=\overline{\nabla}_{X_i}d\varphi(X_j)-d\varphi(\nabla_{X_i}X_j)\in 
 T^{\perp}M,
 \\
 &\overline{\nabla}_{X_i}\tau_b(\varphi)\in T^{\bot}M, 
 \end{aligned}
 \right.
 \end{equation}
 where $\nabla$ is the Tanaka-Webster connection 
 and $\nabla_{X_i}X_j\in H(M)$. 
 Then, we have, in the first term of the second sum of (5.7),
 \begin{align}
 \sum_{i,j=1}^{2n}
 \langle
 \overline{\nabla}_{X_i}&\tau_b(\varphi),
 \overline{\nabla}_{X_i}d\varphi(X_j)\rangle\,
 d\varphi(X_j)\nonumber\\
 &=
 \sum_{i,j=1}^{2n}\left\langle 
 \overline{\nabla}_{X_i}\tau_b(\varphi), 
 (\widetilde{\nabla}_{X_i}d\varphi)(X_j)
 +d\varphi(\nabla_{X_i}X_j)\right\rangle\, d\varphi(X_j)\nonumber\\
 &=\sum_{i,j=1}^{2n}
 \langle\overline{\nabla}_{X_i}\tau_b(\varphi),d\varphi(\nabla_{X_i}X_j)\rangle\,d\varphi(X_j)\nonumber\\
 &=
 \sum_{i,j=1}^{2n}\left\langle\overline{\nabla}_{X_i}\tau_b(\varphi),d\varphi\left(\sum_{k=1}^{2n}\langle \nabla_{X_i}X_j,X_k\rangle \,X_k\right)\right\rangle\,d\varphi(X_j), 
 \end{align}
 because of $\nabla_{X_i}X_j\in H(M)$. 
 Since the Tanaka-Webster connection $\nabla$ satisfies 
 $\nabla g_{\theta}=0$, we have 
 \begin{align}
 \langle \nabla_{X_i}X_j,X_k\rangle
 =X_i\,\langle X_j,X_k\rangle-\langle X_j,\nabla_{X_i}X_k\rangle=-\langle X_j,\nabla _{X_i}X_k\rangle.\nonumber
 \end{align}
 Thus, (5.9) turns to 
 \begin{align}
 \sum_{i,j=1}^{2n}
 \langle \overline{\nabla}_{X_i}\tau_b(\varphi),&\overline{\nabla}_{X_i}d\varphi(X_j)\rangle\,d\varphi(X_j)\nonumber\\
 &=
 \sum_{i,j,k=1}^{2n} 
 \langle\nabla_{X_i}X_j,X_k\rangle\,
 \langle\overline{\nabla}_{X_i}\tau_b(\varphi),d\varphi(X_k)\rangle\,d\varphi(X_j)\nonumber\\
 &=-\sum_{i,j,k=1}^{2n}\langle X_j,\nabla_{X_i}X_k\rangle\,\,\langle \overline{\nabla}_{X_i}\tau_b(\varphi),d\varphi(X_k)\rangle\,d\varphi(X_j)\nonumber\\
 &=-\sum_{i,k=1}^{2n}\langle\overline{\nabla}_{X_i}\tau_b(\varphi),d\varphi(X_k)\rangle\,d\varphi(\nabla_{X_i}X_k).
 \end{align}
 \par
 (The forth step) \quad 
 By inserting (5.10) into (5.7), the second sum of 
 (5.7) turns to 
 \begin{align}
 &\sum_{i,j=1}^{2n} 
 \left[
 \begin{aligned}
 &\langle\overline{\nabla}_{X_i}\tau_b(\varphi),\overline{\nabla}_{X_i}
(d\varphi(X_j))\rangle\,d\varphi(X_j)\\
&+\langle\overline{\nabla}_{X_i}\tau_b(\varphi),d\varphi(X_j)\rangle\,\overline{\nabla}_{X_i}d\varphi(X_j) 
\end{aligned}
 \right]\nonumber
 \\
 &=
 \sum_{i,j=1}^{2n}\left[
 \begin{aligned}
 &\langle \overline{\nabla}_{X_i}\tau_b(\varphi),d\varphi(X_j)\rangle\,\overline{\nabla}_{X_i}d\varphi(X_j)\nonumber\\
 &-\langle\overline{\nabla}_{X_i}\tau_b(\varphi),d\varphi(X_j)\rangle\,d\varphi(\nabla_{X_i}X_j)
 \end{aligned}\nonumber
 \right]\nonumber\\
 &=
 \sum_{i,j=1}^{2n}\langle\overline{\nabla}_{X_i}\tau_b(\varphi),d\varphi(X_j)\rangle\,\left(
 \widetilde{\nabla}_{X_i}d\varphi\right)(X_j).
 \end{align}
 Thus, by (5.7) and (5.11), we obtain Lemma 5.3. \qed
 \vskip0.6cm\par
 ({\em Proof of Lemma 5.4}) \quad 
 We will calculate the right hand side of (5.2) in Lemma 5.3. 
 The proof is divided into several steps. 
 \medskip\par
 (The first step)\quad We first note that 
 \begin{align}
 \langle \overline{\nabla}_{X_i}\tau_b(\varphi),d\varphi(X_j)\rangle+
 \langle
 \tau_b(\varphi),\overline{\nabla}_{X_i}d\varphi(X_j)\rangle &=X_i\,\langle \tau_b(\varphi),d\varphi(X_j)\rangle\nonumber\\
 &=0.
 \end{align}
 Thus, by (5.12), 
 we have 
 \begin{align}
 \langle\overline{\nabla}_{X_i}\tau_b(\varphi),d\varphi(X_j)\rangle&=-\langle\tau_b(\varphi),\overline{\nabla}_{X_i}d\varphi(X_j)\rangle\nonumber\\
 &=-\langle\tau_b(\varphi),\overline{\nabla}_{X_i}d\varphi(X_j)-d\varphi(\nabla_{X_i}X_j\rangle
 \nonumber\\
 &=
 -\langle\tau_b(\varphi),\overline{\nabla}_{X_i}d\varphi(X_j)-d\varphi(\nabla^{g_{\theta}}{}_{X_i}X_j)\rangle \nonumber\\
 &\qquad\qquad\qquad\qquad 
 \,\, (\mbox{by}\,\langle\tau_b(\varphi),d\varphi(T)\rangle=0)
 \nonumber\\
 &=-\langle\tau_b(\varphi),\left(
 \widetilde{\nabla}_{X_i}d\varphi\right)(X_j)\rangle.
 \end{align}
 \medskip\par
 (The second step)\quad 
 By differentiating (5.12), we have 
 \begin{align}
 \left\langle 
 \overline{\nabla}_{X_i}
 \left(
 \overline{\nabla}_{X_i}\tau_b(\varphi)
 \right), 
 d\varphi(X_j)
 \right\rangle
 &+
 2\left\langle
 \overline{\nabla}_{X_i}\tau_b(\varphi),
 \overline{\nabla}_{X_i}d\varphi(X_j)
 \right\rangle\nonumber\\
 &+\left\langle
 \tau_b(\varphi),
 \overline{\nabla}_{X_i}
 \left(
 \overline{\nabla}_{X_i}d\varphi(X_j)
 \right)
 \right\rangle=0. 
 \end{align} 
 And we have 
 \begin{align}
 \left\langle
 \overline{\nabla}_{\nabla_{X_i}X_i}\tau_b(\varphi),d\varphi(X_j)
 \right\rangle
&+
 \left\langle\tau_b(\varphi),
 \overline{\nabla}_{\nabla_{X_i}X_i}d\varphi(X_j)\right\rangle\nonumber\\
 &=\nabla_{X_i}X_i\,\langle\tau_b(\varphi),d\varphi(X_j)\rangle=0. 
 \end{align}
 Thus, by (5.14) and (5.15), we have 
 \begin{align}
 \langle-\Delta_b\,(\tau_b(\varphi)),&d\varphi(X_j)\rangle
 +2\langle\overline{\nabla}_{X_i}\tau_b(\varphi),\overline{\nabla}_{X_i}d\varphi(X_j)\rangle\nonumber\\
 &+\langle\tau_b(\varphi),-\Delta_b(d\varphi(X_j))\rangle=0.
 \end{align}
 \medskip\par
 (The third step)\quad 
 For the second term of the left hand side of (5.16), 
 we have 
 \begin{align}
 2\langle\overline{\nabla}_{X_i}\tau_b(\varphi),\overline{\nabla}_{X_i}d\varphi(X_j)\rangle
 =-2\left\langle \tau_b(\varphi),\left(\widetilde{\nabla}_{X_i}d\varphi\right)(\nabla_{X_i}X_j)\right\rangle. 
 \end{align}
 Because, the left hand side of (5.17) is 
 \begin{align}
 2\langle
 \overline{\nabla}_{X_i}\tau_b(\varphi),\overline{\nabla}_{X_i}d\varphi(X_j)\rangle&=
 2\left\langle
 \overline{\nabla}_{X_i}\tau_b(\varphi),
 \left(\widetilde{\nabla}_{X_i}d\varphi
 \right)(X_j)+d\varphi
 \left(
 \nabla^{g_{\theta}}_{X_i}X_j
 \right)
 \right\rangle\nonumber\\
 &=
 2\left\langle
 \overline{\nabla}_{X_i}\tau_b(\varphi),
 \left(\widetilde{\nabla}_{X_i}d\varphi
 \right)(X_j)+d\varphi
 \left(
 \nabla_{X_i}X_j
 \right)
 \right\rangle\nonumber\\
 &=2\left\langle
 \overline{\nabla}_{X_i}\tau_b(\varphi),
d\varphi
 \left(
 \nabla_{X_i}X_j
 \right)
 \right\rangle\nonumber\\
 &=-2
 \left\langle \tau_b(\varphi),\overline{\nabla}_{X_i}d\varphi(\nabla_{X_i}X_j)\right\rangle
\nonumber\\
&\qquad\qquad\qquad\qquad (\mbox{by $\langle\tau_b(\varphi),d\varphi(\nabla_{X_i}X_j)\rangle=0$})
 \nonumber\\
 &=-2\left\langle\tau_b(\varphi),
 \left(\widetilde{\nabla}_{X_i}d\varphi\right)(\nabla_{X_i}X_j)\right\rangle
 \end{align} 
 which is the right hand side of (5.17). 
 In the last step of (5.18), we used the equality
 $\langle\tau_b(\varphi),\overline{\nabla}_{X_i}d\varphi(T)\rangle=0$ 
 which follows from that $\overline{\nabla}_{X_i}d\varphi(T)$ is tangential. 
  \medskip\par
  (The forth step)\quad 
  For the third term of the left hand side of (5.16), we have 
  \begin{equation}
  \langle \tau_b(\varphi),-\Delta_b(d\varphi(X_j))\rangle=
  \big\langle\tau_b(\varphi),\left(-\widetilde{\Delta}_bd\varphi\right)(X_j)
  +2
  \sum_{k=1}^{2n}(\widetilde{\nabla}_{X_k}d\varphi)(\nabla_{X_k} X_j  )\big\rangle.
  \end{equation} 
  Because, by the definition of $\Delta_b$, we have 
   \begin{align}
   \langle 
   \tau_b(\varphi),
   &-\Delta_b(d\varphi(X_j))\rangle \nonumber\\
   &=
   \langle\tau_b(\varphi), 
   \sum_{k=1}^{2n}\left\{
   \overline{\nabla}_{X_k}\left(
   \overline{\nabla}_{X_k}d\varphi(X_j)
   \right)
   -\overline{\nabla}_{\nabla_{X_k}X_k}d\varphi(X_j)
   \right\}\rangle\nonumber\\
   &=\big\langle 
   \tau_b(\varphi),
   \sum_{k=1}^{2n}\big\{
   \overline{\nabla}_{X_k}\big(
   (\widetilde{\nabla}_{X_k}d\varphi)(X_j)+d\varphi(\nabla_{X_k}X_j)\big)\nonumber\\
   &\qquad\qquad\qquad
   -\big(\widetilde{\nabla}_{\nabla_{X_k}X_k}d\varphi\big)(X_j)-d\varphi\big(
   \nabla_{\nabla_{X_k}X_k}X_j
   \big)
   \big\}\big\rangle\nonumber\\
   &=
   \big\langle 
   \tau_b(\varphi),\sum_{k=1}^{2n}\big\{
   \big(
   \widetilde{\nabla}_{X_k}\widetilde{\nabla}_{X_k}d\varphi
   \big)(X_j)+
   \big(\widetilde{\nabla}_{X_k}d\varphi\big)(\nabla_{X_k}X_j)\nonumber\\
   &\qquad \qquad\qquad
   +\big(\widetilde{\nabla}_{X_k}d\varphi\big)(\nabla_{X_k}X_j)+d\varphi(\nabla_{X_k}\nabla_{X_k}X_j)\nonumber\\
   &\qquad\qquad\qquad
   -\big(
   \widetilde{\nabla}_{\nabla_{X_k}X_k}
   d\varphi\big)(X_j)
   \big\rangle\big\}
   \nonumber\\
   &=
   \big\langle\tau_b(\varphi),\big(-\widetilde{\Delta}_bd\varphi\big)(X_j)
   +2\sum_{k=1}^{2n}\big(\widetilde{\nabla}_{X_k}d\varphi
   \big)(\nabla_{X_k}X_j)\big\rangle,\nonumber
   \end{align}
   which is (5.19). To get the last equality of the above, 
   we used the following equations: \quad for all $X\in H(M)$, it holds that 
   \begin{equation}
   \langle\tau_b(\varphi),(\overline{\nabla}_{X}d\varphi)(T)\rangle=
   \langle\tau_b(\varphi),d\varphi(X)\rangle=\langle\tau_b(\varphi),(\widetilde{\nabla}_{X}d\varphi)(T)\rangle=0.
   \end{equation}
  To get (5.20), 
  due to the admissibility of $\varphi$, 
  we have 
  $$(\widetilde{\nabla}_Xd\varphi)(T)=(\overline{\nabla}_Xd\varphi)(T)-d\varphi(\nabla^{g_{\theta}}_XT)=B_{\varphi}(X,T)=0,$$ 
  and then, 
  $(\overline{\nabla}_Xd\varphi)(T)$ is tangential 
  for all $X\in H(M)$.  We have (5.20), and then (5.19). 
  \medskip\par
  (The fifth step)\quad 
 Then, the right hand side of (5.19) is equal to 
  \begin{align}
  \langle\tau_b(\varphi),&(-\widetilde{\Delta}_bd\varphi)(X_j)+2\sum_{k=1}^{2n}(\widetilde{\nabla}_{X_k}d\varphi)(\nabla_{X_k}X_j)\rangle\nonumber\\
  &=
  \langle\tau_b(\varphi),\overline{\nabla}_{X_j}\tau(\varphi)\nonumber\\
  &\qquad 
  -\sum_{k=1}^{2n}\left\{
  R^h(d\varphi(X_j),d\varphi(X_k))d\varphi(X_k)-
  d\varphi(R^{g_{\theta}}(X_j,X_k)X_k)\right\}\nonumber\\
  &\qquad 
  -\sum_{k=1}^{2n}\left\{
  R^h(d\varphi(X_j),d\varphi(T))d\varphi(T)-d\varphi((R^{g_{\theta}}(X_j,T)T
  \right\}\nonumber\\
  &\qquad-\widetilde{\nabla}_T\,\widetilde{\nabla}_T\,d\varphi(X_j)\nonumber\\
  &\qquad+2(\widetilde{\nabla}_{X_k}d\varphi)(\nabla_{X_k}X_j)\rangle, 
  \end{align}
  which follows from the formula (3.11).  
  \par
  Here, notice that 
  \begin{align}
  \langle\tau_b(\varphi),\overline{\nabla}_{X_j}\tau(\varphi)\rangle=0.
  \end{align}
  Because $\tau_b(\varphi)$ is normal, and 
  $\overline{\nabla}_{X_j}\tau(\varphi)$ is tangential. 
   And also we have 
   \begin{align}
   &\langle \tau_b(\varphi),\widetilde{\nabla}_T(d\varphi(X_j)\rangle=0,\\
   &\langle\tau_b(\varphi),\widetilde{\nabla}_T\,\widetilde{\nabla}_T\,d\varphi(X_j)\rangle=0.
   \end{align}
   To see (5.23), since we assume $\varphi$ is an admissible isometric immersion, 
   we have 
   \begin{align}
   \widetilde{\nabla}_T\,d\varphi(X_j)=\nabla^h_TX_j=\nabla^{g_{\theta}}_TX_j+B_{\varphi}(T,X_j)=\nabla^{g_{\theta}}_TX_j
   \end{align}
   \par\noindent 
   which is tangential, so that we have (5.23). 
   Furthermore, to see (5.24), we have 
   \begin{align}
   \widetilde{\nabla}_T\,\widetilde{\nabla}_T\,d\varphi(X_j)&=\widetilde{\nabla}(\nabla^{g_{\theta}}_TX_j)\nonumber\\
   &=\nabla^h_T(\nabla^{g_{\theta}}_TX_j)\nonumber\\
   &=\nabla^{g_{\theta}}_T(\nabla^{g_{\theta}}_TX_j)+B(T,\nabla^{g_{\theta}}_TX_j).
   \end{align}
   Here, 
   for every $X\in H(M)$, 
   $$\nabla^{g_{\theta}}_TX\in H(M).$$
   \par\noindent
   Indeed, since 
   $g_{\theta}(T,X)=0$, and $\nabla^{g_{\theta}}_TT=0$ (cf. \cite{DT}, pp. 47, and 48), 
   \begin{align*}
   g_{\theta}(T,\nabla^{g_{\theta}}_TX)&=T(g_{\theta}(T,X))-g_{\theta}(\nabla^{g_{\theta}}_TT,X)=0,
   \end{align*}
   which implies $\nabla^{g_{\theta}}_TX\in H(M)$. 
   Thus, the admissibility implies that the second term of (5.26) vanishes. 
Thus, the right hand side of (5.26) is tangential, which implies (5.24).  
\par
Therefore, we obtain 
\begin{align}
  \langle\tau_b(\varphi),&(-\widetilde{\Delta}_bd\varphi)(X_j)+2\sum_{k=1}^{2n}(\widetilde{\nabla}_{X_k}d\varphi)(\nabla_{X_k}X_j)\rangle\nonumber\\
  &=
-\sum_{k=1}^{2n}\langle
\tau_b(\varphi),R^h(d\varphi(X_j),d\varphi(X_k))d\varphi(X_k)\nonumber\\
&\quad-\sum_{k=1}^{2n}\langle\tau_b(\varphi),
R^h(d\varphi(X_j),d\varphi(T))d\varphi(T)\rangle\nonumber\\
&\quad +2\sum_{k=1}^{2n}\langle\tau_b(\varphi),(\widetilde{\nabla}_{X_k}d\varphi)(\nabla_{X_k}X_j)\rangle.
\end{align}
\par
(The sixth step)\quad Now, return to (5.16), by using (5.17), (5.19) and (5.27), we have 
\begin{align}
0&= \langle-\Delta_b\,(\tau_b(\varphi)),
d\varphi(X_j)\rangle
 +2\langle\overline{\nabla}_{X_i}\tau_b(\varphi),
 \overline{\nabla}_{X_i}d\varphi(X_j)\rangle\nonumber
 \\
 &\qquad 
 +\langle\tau_b(\varphi),-\Delta_b(d\varphi(X_j))\rangle\nonumber\\
 &=\langle-\Delta_b(\tau_b(\varphi)),d\varphi(X_j)\rangle-2\langle\tau_b(\varphi),(\widetilde{\nabla}_{X_i}d\varphi)(\nabla_{X_i}X_j)\rangle\nonumber\\
& \qquad
 +\langle\tau_b(\varphi),-\Delta_b(d\varphi(X_j))\rangle
 \nonumber\\
 &=\langle-\Delta_b(\tau_b(\varphi)),d\varphi(X_j)\rangle
 -2\langle\tau_b(\varphi),(\widetilde{\nabla}_{X_i}d\varphi)(\nabla_{X_i}X_j)\rangle
 \nonumber\\
 &\qquad 
 +\langle\tau_b(\varphi),(-\widetilde{\Delta}_bd\varphi)(X_j)+2\sum_{k=1}^{2n}(\widetilde{\nabla}_{X_k}d\varphi)(\nabla_{X_k}X_j)\rangle
 \nonumber\\
 &=
 \langle-\Delta_b\,(\tau_b(\varphi)),
d\varphi(X_j)\rangle
 -2\langle
\tau_b(\varphi),
 (\widetilde{\nabla}_{X_i}
 d\varphi)
 (\nabla_{X_i}X_j)\rangle
 \nonumber\\
 &\qquad 
 +\bigg\langle\tau_b(\varphi),-\sum_{k=1}^{2n}R^h(d\varphi(X_j),d\varphi(X_k))d\varphi(X_k)
 \nonumber\\
 &\qquad\qquad \qquad
 -R^h(d\varphi(X_j),d\varphi(T))d\varphi(T)\nonumber\\
 &\qquad\qquad \qquad
 +2\sum_{k=1}^{2n}(\widetilde{\nabla}_{X_k}d\varphi)(\nabla_{X_k}X_j)\bigg\rangle\nonumber\\
 &=
 \langle-\Delta_b\,(\tau_b(\varphi)),
d\varphi(X_j)\rangle\nonumber\\
&\qquad 
+\bigg\langle\tau_b(\varphi),
-\sum_{k=1}^{2n}R^h(d\varphi(X_j),d\varphi(X_k))d\varphi(X_k)
 \nonumber\\
 &\qquad\qquad \qquad
 -R^h(d\varphi(X_j),d\varphi(T))d\varphi(T)\bigg\rangle. 
\end{align}
\par
(The seventh step) \quad 
Inserting (5.28) into (5.2) of Lemma 5.3, 
we obtain 
\begin{align}
-\Delta_b(\tau_b(\varphi))&=
\langle\tau_b(\varphi),\sum_{k=1}^{2n}R^h(d\varphi(X_j),d\varphi(X_k))d\varphi(X_k)\rangle\,d\varphi(X_j)\nonumber\\
&\quad +\langle \tau_b(\varphi),R^h(d\varphi(X_j),d\varphi(T))d\varphi(T)\rangle\,d\varphi(X_j)\nonumber\\
&\quad +\langle\overline{\nabla}_{X_i}\tau_b(\varphi),d\varphi(X_j)\rangle\,(\widetilde{\nabla}_{X_i}d\varphi)(X_j). 
\end{align}
\par
At last, for the third term of (5.29), 
we have 
\begin{align}
\langle\overline{\nabla}_{X_i}\tau_b(\varphi),d\varphi(X_j)\rangle&=
X_i\,\langle\tau_b(\varphi),d\varphi(X_j)\rangle-\langle\tau_b(\varphi),\overline{\nabla}_{X_i}d\varphi(X_j)\rangle\nonumber\\
&=-\langle\tau_b(\varphi),(\widetilde{\nabla}_{X_i}d\varphi)(X_j)\rangle.
\end{align}
Together with (5.29) and (5.30), we have (5.3) of Lemma 5.4. 
\qed
\bigskip\par
Due to Lemma 5.4 and the definition of biharmonicity, we obtain immediately 
\begin{th} 
Let 
$\varphi$ be an admissible isometric immersion of a strongly pseudoconvex $CR$ manifold 
$(M,g_{\theta})$ into another Riemannian manifold $(N,h)$ 
whose pseudo mean curvature vector field along $\varphi$ is parallel. Then, 
$\varphi$ is pseudo biharmonic if and only if 
\begin{align}
\tau_{b,2}(\varphi):=\Delta_b\big(\tau_b(\varphi)\big)-\sum_{j=1}^{2n}R^h\big(\tau_b(\varphi),d\varphi(X_j)\big)d\varphi(X_j)=0
\end{align} 
if and only if 
\begin{align}
&-\sum_{j,k=1}^{2n}h\big(\tau_b(\varphi),R^h(d\varphi(X_j),d\varphi(X_k))d\varphi(X_k)\big)\,d\varphi(X_j)\nonumber\\
&-\sum_{j=1}^{2n}h\big(\tau_b(\varphi),R^h(d\varphi(X_j),d\varphi(T))d\varphi(T)\big)\,d\varphi(X_j)\nonumber\\
&+\sum_{j,k=1}^{2n}h\big(\tau_b(\varphi),B_{\varphi}(X_j,X_k)\big)\,B_{\varphi}(X_j,X_k)\nonumber\\
&-\sum_{j=1}^{2m}R^h\big(\tau_b(\varphi),d\varphi(X_j)\big)d\varphi(X_j)=0,
\end{align}
where $\{X_j\}_{j=1}^{2n}$ is an orthonormal frame field 
of $(H(M), g_{\theta})$. 
\end{th}
\medskip\par
\begin{rem}
Due to \cite{J}, p. 220, Lemma 10,  and the definition of 
bi-tension field $\tau_2(\varphi)$ for an isometric immersion $\varphi$ of a Riemannian manifold $(M,g)$ into another Riemannian manifold $(N,h)$, 
we can also obtain immediately the following useful theorem: 
\par
{\bf Theorem}\quad 
Let $\varphi$ be an isometric immersion of a Riemannian manifold $(M^m,g)$ into another Riemannian manifold $(N^n,h)$ 
whose mean curvature vector field along $\varphi$ is parallel. Let $\{e_j\}_{j=1}^m$ be an orthonormal frame field of $(M,g)$. 
Then, $\varphi$ is biharmonic if and only if 
\begin{align}
\tau_2(\varphi):=\overline{\Delta}\big(\tau(\varphi)\big)-\sum_{j=1}^mR^h\big(\tau(\varphi),e_j\big)e_j=0
\end{align}  
if and only if 
\begin{align}
&-\sum_{j,k=1}^mh\big(\tau(\varphi),R^h(d\varphi(e_j),d\varphi(e_k))d\varphi(e_k)\big)\,d\varphi(e_j)\nonumber\\
&+\sum_{j,k=1}^mh\big(\tau(\varphi),B_{\varphi}(e_j,e_k)\big)\,B_{\varphi}(e_j,e_k)\nonumber\\
&-\sum_{j=1}^mR^h\big(\tau(\varphi),d\varphi(e_j)\big)d\varphi(e_j)=0.
\end{align}
\end{rem}
\medskip\par
\section{Isometric immersions into the unit sphere.}
In this section, we treat with 
admissible isometric immersions 
of $(M^{2n+1},g_{\theta})$ into the unit sphere 
$(N,h)=S^{2n+2}(1)$ 
with parallel pseudo mean curvature vector field 
with non-zero pseudo mean curvature. 
 \par
 The curvature tensor field $R^h$ of the target space 
 $(N,h)=S^{2n+2}(1)$ satifies that 
 \begin{align}
 R^h(X,Y)Z=h(Z,Y)\,X-h(Z,X)\,Y
 \end{align}
 for all vector fields $X,\,Y,\,Z$ on $N$. Then, we have 
 \begin{align}
 &\left(
 R^h(d\varphi(X_j),d\varphi(X_k))\,d\varphi(X_k)
 \right)^{\perp}=0,\\
 &
 \left(R^h(d\varphi(X_j),d\varphi(T))d\varphi(T)
 \right)^{\perp}=0, 
 \end{align}
 for all $i,j=1,\,\cdots,2n$. 
 Therefore, we obtain by (5.3) in Lemma 5.4, 
 \begin{align}
 -\Delta_b (\tau_b(\varphi))
 =-\left\langle\tau_b(\varphi),\left(
 \widetilde{\nabla}_{X_i}d\varphi
 \right)(X_j)\right\rangle\,\left(\widetilde{\nabla}_{X_i}d\varphi\right)(X_j).
 \end{align}
 \par
 On the other hand, 
 we have 
 \begin{align}
 \sum_{k=1}^{2n}&R^h(\tau_b(\varphi),d\varphi(X_k))d\varphi(X_k)\nonumber\\
 &=\sum_{k=1}^{2n}h(d\varphi(X_k),d\varphi(X_k))\,
\tau_b(\varphi)
 -\sum_{k=1}^{2n}h(d\varphi(X_k),\tau_b(\varphi))\,\tau_b(\varphi)\nonumber\\
 &=2n\,\tau_b(\varphi).
 \end{align}
 \par
 Now, let us recall the pseudo biharmonicity of $\varphi$ is equivalent to 
 that 
 \begin{align}
 -\Delta_b(\tau_b(\varphi))+\sum_{k=1}^{2n}R^h(\tau_b(\varphi),d\varphi(X_k))\,d\varphi(X_k)=0
 \end{align}
 which is equivalent to that 
 \begin{align}
 -\left\langle \tau_b(\varphi), 
 \left(\widetilde{\nabla}_{X_i}d\varphi\right)(X_j)\right\rangle\,
 \left(
 \widetilde{\nabla}_{X_i}d\varphi\right)(X_j)+2n\,\tau_b(\varphi)=0.
 \end{align}
 \par
 For our immersion $\varphi:\,(M,g_{\theta})\rightarrow S^{2n+2}(1)$, 
 let $\xi$ be the unit normal vector filed on $M$ along $\varphi$, 
 we have by definition, 
 \begin{align}
 \left(
 \widetilde{\nabla}_{X_i}d\varphi\right)(X_j)=B_{\varphi}(X_i,X_j)=H_{ij}\,\xi.
 \end{align}
 Then, we have by definition of $\tau_b(\varphi)$, 
 \begin{align}
 \tau_b(\varphi)=\sum_{i=1}^{2n}\left(\widetilde{\nabla}_{X_i}d\varphi\right)(X_i)
 =\left(\sum_{i=1}^{2n}H_{ii}\right)\,\xi.
 \end{align}
 Therefore, we have 
 \begin{align}
 \Vert\tau_b(\varphi)\Vert^2=\left(\sum_{i=1}^{2n}H_{ii}\right)^2\,\Vert\xi\Vert^2=
 \left(\sum_{i=1}^{2n}H_{ii}\right)^2.
 \end{align}
 By the admissibility, we have 
 \begin{align}
 \Vert B_{\varphi}\Vert^2&=\sum_{i,j=1}^{2n}\Vert B_{\varphi}(X_i,X_j)\Vert^2
 +2\sum_{i=1}\Vert B_{\varphi}(X_i,T)\Vert^2+\Vert B_{\varphi}(T,T)\Vert^2\nonumber\\
 &=\sum_{i,j=1}^{2n}\Vert H_{ij}\,\xi\Vert^2+\Vert B_{\varphi}(T,T)\Vert^2\nonumber\\
 &=\sum_{i,j=1}^{2n}H_{ij}{}^2+\Vert B_{\varphi}(T,T)\Vert^2. 
 \end{align} 
 Due to (6.7), (6.8) and (6.9), the biharmonicity of $\varphi$ is equivalent to that 
 \begin{align}
 0=-\left\langle\left(
 \sum_{k=1}^{2n}H_{kk}\right)\,\xi,H_{ij}\,\xi\right\rangle\,H_{ij}\,\xi
 +2n\,\left(\sum_{k=1}^{2n}H_{kk}\right)\,\xi
 \end{align}
 which is equivalent to that 
 \begin{align}
 0&=\left(\sum_{k=1}^{2n}H_{kk}
 \right)\,\left\{-\sum_{i,j=1}^{2n}H_{ij}{}^2+2n
 \right\}\nonumber\\
 &=\Vert\tau_b(\varphi)\Vert\,\left\{
 -\Vert B_{\varphi}\Vert^2+\Vert B_{\varphi}(T,T)\Vert^2+2n
 \right\}
 \end{align}
 by (6.11). By our assumption of non-zero pseudo mean curvature, 
 $\Vert \tau_b(\varphi)\Vert \not=0$ at every point,  we obtain the following equivalence relation: 
 $\varphi$ is pseudo biharmonic if and only if 
 \begin{align}
 \Vert B_{\varphi}\Vert^2=
 \Vert B_{\varphi}(T,T)\Vert^2+2n
 \end{align}
 at every point in $M$. 
 \par
 By summing up the above, we obtain the following theorem:
 \begin{th}
 Let $\varphi$ be an sdmissible isometric immersion of 
 a strictly pseudoconvex $CR$ manifold $(M,g_{\theta})$ 
 into the unit sphere $(N,h)=S^{2n+2}(1)$ Assume that  the pseudo mean curvature vector field is parallel with non-zero pseudo mean curvature. 
 Then, $\varphi$ is pseudo biharmonic if and only if 
 \begin{align}
 \Vert B_{\varphi}\Vert^2=\Vert B_{\varphi}(T,T)\Vert^2+2n.
 \end{align}
 \end{th}
 \bigskip\par
The admissibility condition is that: 
$d\varphi(T)$ is the principal curvature 
vector field along $\varphi$ with some principal curvature, say $\lambda_{2n+1}$. I.e., 
$$
A_{\xi}T=\lambda_{2n+1}\,T.
$$ 
Then, we have 
$$
\Vert B_{\varphi}\Vert^2=\sum_{i=1}^{2n+1}\lambda_i{}^2,\quad 
\mbox{and}\,\,\Vert B_{\varphi}(T,T)\Vert^2=\lambda_{2n+1}.
$$
By Theorem 6.1, we have immediately 
Thus, we obtain 
\medskip\par
\begin{cor}
Let $\varphi:\,(M^{2n+1},g_{\theta})\rightarrow S^{2n+2}(1)$ be an isometric immersion whose the pseudo mean curvature vector field is parallel and 
has non-zero pseudo mean curvature. Then, 
$\varphi$ is pseudo biharmonic if and only if it holds that 
\begin{equation}
\sum_{i=1}^{2n}\lambda_i{}^2=2n
\end{equation}
which is equivalent to that 
\begin{align}
\bigg\Vert B_{\varphi}\bigg\vert_{H(M)\times H(M)}\bigg\Vert^2=2n,
\end{align}
where 
$B_{\varphi}\vert_{H(M)\times H(M)}$ is the restriction of $B_{\varphi}$ to the subspace $H(M)$ of the tangent space $T_xM$ $(x\in M)$. 
\end{cor}
\bigskip\par
\section{Isometric immersions to the complex projective space.}
In this section, we will consider 
admissible isometric immersions 
of $(M^{2n+1},g_{\theta})$ into 
the complex projective space 
$(N,h)={\mathbb P}^{n+1}(c)$ ($c>0$) 
whose 
mean curvature vector field is parallel with non-zero pseudo mean curvature. 
\medskip\par
{\bf 7.1.} \quad 
Let us recall that the curvature tensor field $(N,h)={\mathbb P}^{n+1}(c)$ is 
given by 
\begin{align}
R^h(U,V)W=\frac{c}{4}\,\big\{
&h(V,W)\,U-h(U,W)\,V\nonumber\\
&+h(JV,W)\,JU-h(JU,W)\,JV+2h(U,JV)\,JW\big\},
\end{align}
where 
$J$ is the adapted complex structure of ${\mathbb P}^{n+1}(c)$, and 
$U,\,V$ and $W$ are vector fields on ${\mathbb P}^{n+1}(c)$, respectively. 
Therefore, we have 
\begin{align}
&R^h(d\varphi(X_j),d\varphi(X_k))d\varphi(X_k)\nonumber\\
&=\frac{c}{4}\big\{
h(d\varphi(X_k),d\varphi(X_k))\,d\varphi(X_j)-h(d\varphi(X_j),d\varphi(X_k))\,d\varphi(X_k)\nonumber\\
&\quad +h(J\,d\varphi(X_k),d\varphi(X_k))\,Jd\varphi(X_j)-h(Jd\varphi(X_j),d\varphi(X_k))\,Jd\varphi(X_k)\nonumber\\
&\quad +2h(d\varphi(X_j),Jd\varphi(X_k))\,Jd\varphi(X_k)\big\}\nonumber\\
&=\frac{c}{4}\bigg\{
d\varphi(X_j)-\delta_{jk}\,d\varphi(X_k)+3\,h(d\varphi(X_j),Jd\varphi(X_k))\,Jd\varphi(X_k)\bigg\}.
\end{align}
We show first 
\begin{align}
\sum_{j,k=1}^{2n}
&h\left(
\tau_b(\varphi),R^h(d\varphi(X_j),d\varphi(X_k))\,
d\varphi(X_k)\right)\,d\varphi(X_j)\nonumber\\
&=-\frac{3c}{4}\,h(\tau_b(\varphi),J\,d\varphi(T))\,
\big(
J\,d\varphi(T)
\big)^{\top}
\nonumber\\
&\quad-\frac{3c}{4}\sum_{j=1}^{2n}h\big(d\varphi(X_j),J\,\big(J\,\tau_b(\varphi)\big)^{\top}\,\big)\,d\varphi(X_j).
\end{align}
Recall here that the tangential part of $Z\in T_{\varphi(x)}N\,\,(x\in M)$ is given by
\begin{align}
Z^{\top}=\sum_{i=1}^{2n}h(Z,d\varphi(X_i))\,d\varphi(X_i)+h(Z,d\varphi(T))\,d\varphi(T).
\end{align} 
Since 
$h(\tau_b(\varphi),d\varphi(X_j))=0$ 
$(j=1,\cdots,2n)$, and (7.2), 
one can calculate the left hand side of (7.3) 
as follows: 
\begin{align}
\sum_{j,k=1}^{2n}&h\big(\tau_b(\varphi),R^h(d\varphi(X_j),d\varphi(X_k))d\varphi(X_k)\big)\,d\varphi(X_j)\nonumber\\
&=\frac{3c}{4}\,\sum_{j,k=1}^{2n}h(\tau_b(\varphi), 
h(d\varphi(X_j),J\,d\varphi(X_k))\,h(\tau_b(\varphi),J\,d\varphi(X_k))\,d\varphi(X_j)\nonumber\\
&=\frac{3c}{4}\,\sum_{j,k=1}^{2n}h(J\,d\varphi(X_j),d\varphi(X_k))\,h(J\,\tau_b(\varphi),d\varphi(X_k))\,d\varphi(X_j)\nonumber\\
&=\frac{3c}{4}\,\sum_{j=1}^{2n}h\big(
J\,d\varphi(X_j),\sum_{k=1}^{2n}h(J\,\tau_b(\varphi),
d\varphi(X_k))\,d\varphi(X_k))\,d\varphi(X_j)
\nonumber\\
&=
\frac{3c}{4}
\sum_{j=1}^{2n}
h\big(J\,d\varphi(X_j),
(\,J\tau_b(\varphi))^{\top}
-h(J\,\tau_b(\varphi),d\varphi(T))\,d\varphi(T)
\big)\,d\varphi(X_j)\nonumber\\
&=\frac{3c}{4}\sum_{j=1}^{2n}h(J\,d\varphi(X_j),\,(J\,\tau_b(\varphi)\,)^{\top}\,)\,d\varphi(X_j)\nonumber\\
&\quad 
+\frac{3c}{4}\,h(J\tau_b(\varphi),d\varphi(T))\,\sum_{j=1}^{2n}h(d\varphi(X_j),J\,d\varphi(T))\,d\varphi(X_j)\nonumber\\
&=
 -\frac{3c}{4}\,\sum_{j=1}^{2n}h(d\varphi(X_j),J\,(J\tau_b(\varphi)\,)^{\top})\,d\varphi(X_j)
 \nonumber\\
 &\quad
-\frac{3c}{4}\,h(\tau_b(\varphi),J\,d\varphi(T)\,)\,(J\,d\varphi(T))^{\top}.
\end{align}
Then, (7.5) is just  (7.3). 
\par
Second, by a similar way, 
\begin{align}
\sum_{j=1}^{2n}
&\langle\tau_b(\varphi),R^h(d\varphi(X_j),d\varphi(T))\,d\varphi(T)\rangle\,d\varphi(X_j)\nonumber\\
&=\sum_{j=1}^{2n}\langle\tau_b(\varphi), 
\frac{3c}{4}\,h(d\varphi(X_j),J\,d\varphi(T))\,J\,d\varphi(T)\rangle\,d\varphi(X_j)\nonumber\\
&=\frac{3c}{4}\,h(\tau_b(\varphi),J\,d\varphi(T))\,
(J\,d\varphi(T))^{\top} 
\end{align}
in the last equality of (7.6) we used that $h(d\varphi(T),J\,d\varphi(T))=0$. 
\par 
Thus, we have 
\begin{align}
\Delta_b(\tau_b(\varphi))
&=\frac{3c}{4}h(\tau_b(\varphi),Jd\varphi(T))\,(Jd\varphi(T))^{\top}\nonumber\\
&\qquad +\frac{3c}{4}\sum_{j=1}^{2n}h(d\varphi(X_j),
J(J\tau_b(\varphi))^{\top})\,d\varphi(X_j)\nonumber\\
&\qquad 
 -\frac{3c}{4}\,h(\tau_b(\varphi),J\,d\varphi(T))\,(J\,d\varphi(T)\,)^{\top}\nonumber\\
 &\qquad+\langle\tau_b(\varphi),B_{\varphi}(X_i,X_j)\rangle\,B_{\varphi}(X_i,X_j)\nonumber\\
 &=\frac{3c}{4}\sum_{j=1}^{2n}h(d\varphi(X_j),J(\,(J\tau_b(\varphi))^{\top})\,d\varphi(X_j)\nonumber\\
 &\qquad+\langle(\tau_b(\varphi), B_{\varphi}(X_i,X_j)\rangle\,B_{\varphi}(X_i,X_j). 
 \end{align}
 \par
 Therefore, an isometric immersion $\varphi$ is pseudo biharmonic if and only if the pseudo biharmonic map equation folds: 
 \begin{align}
 \Delta_b(\tau_b(\varphi))-\sum_{k=1}^{2n}R^h(\tau_b(\varphi),d\varphi(X_k)\,d\varphi(X_k)=0.
 \end{align}
 By (7.7) and (7.1), (7.8) is equivalent to that the following (7.9) 
 holds: 
 \begin{align}
 &\frac{3c}{4}\sum_{j=1}^{2n}h(d\varphi(X_j),J(\,(J\tau_b(\varphi)^{\top}\,))\,d\varphi(X_j)\nonumber\\
 &\quad+\langle\tau_b(\varphi),B_{\varphi}(X_i,X_j)\rangle\,B_{\varphi}(X_i,X_j)\nonumber\\
 &\quad-\frac{2nc}{4}\,\tau_b(\varphi)+\frac{3c}{4}\sum_{k=1}^{2n}h(d\varphi(X_j),J\tau_b(\varphi))\,J\,d\varphi(X_k)\nonumber\\
 &=0.
 \end{align}
 \medskip\par
 {\bf 7.2.} \quad Let $\xi$ be the unit normal vector field along the admissible isometric immersion 
 $\varphi:\,(M,g_{\theta})\rightarrow {\mathbb P}^{n+1}(c)$ $(c>0)$. 
 \par
  We have immediately 
  \begin{align}
  \left\{
  \begin{aligned}
  &B_{\varphi}(X_i,X_i)=(\widetilde{\nabla}_{X_i}d\varphi)(X_j)=H_{ij}\,\xi,\\
  &\tau_b(\varphi)=\sum_{k=1}^{2n}(\widetilde{\nabla}_{X_k}d\varphi)(X_k)=\left(\sum_{k=1}^{2n}H_{kk}\right)\,\xi,\\
  &J\,\tau_b(\varphi)=\left(\sum_{k=1}^{2n}H_{kk}
  \right)\,J\xi,
  \end{aligned}
  \right.
  \end{align}
  and then, we have 
  \begin{align}
  h(\xi,\,J\,\xi)=0.
  \end{align}
  Indeed, we have 
  \begin{align}
  &-h(J\,\xi,\xi)=h(J\,\xi,J(J\,\xi))=h(\xi,J\,\xi),\nonumber
  \end{align}
 which implies that 
 $h(\xi,J\,\xi)=0$. 
 \par
 Due to (7.11), $J\,\xi$ is tangential.  By (7.10), 
 $J\,\tau_b(\varphi)$ is also tangential.  Therefore, we have 
 \begin{align}
 (\,J\,\tau_b(\varphi)\,)^{\top}=J\,\tau_b(\varphi).
 \end{align}
 In particular, we have 
 \begin{align}
 \sum_{j=1}^{2n}h(d\varphi(X_j),&J\,(J\,\tau_b(\varphi))^{\top}\,)\,)\,d\varphi(X_j)\nonumber\\
 &=
 \sum_{j=1}^{2n}h(d\varphi(X_j),J\,(J\,\tau_b(\varphi))\,)\,d\varphi(X_j)\nonumber\\
 &=-\sum_{j=1}^{2n}h(d\varphi(X_j),\tau_b(\varphi))\,d\varphi(X_j)\nonumber\\
 &=0
 \end{align}
 by using (7.12) and $\tau_b(\varphi)$ is a normal 
 vector field along $\varphi$.    
 \par
 Since $J\,\tau_b(\varphi)$ is tangential, we can write as
 \begin{align}
 J\,\tau_b(\varphi)=
 \sum_{k=1}^{2n}h(d\varphi(X_k),\,J\,\tau_b(\varphi))\,d\varphi(X_k)
+h(d\varphi(T),J\,\tau_b(\varphi))\,d\varphi(T),
\nonumber
 \end{align}
 which implies that 
 \begin{align}
  \sum_{k=1}^{2n}&h(d\varphi(X_k),\,J\,\tau_b(\varphi))\,d\varphi(X_k)\nonumber\\
  &=
 J\,\tau_b(\varphi)
 -h(d\varphi(T),J\,\tau_b(\varphi))\,d\varphi(T). 
 \end{align}
 \par
 Therefore, applying $J$ to (7.14), we have 
 \begin{align}
 \sum_{k=1}^{2n}h(d\varphi(X_k),&J\,\tau_b(\varphi))\,J\,d\varphi(X_k)\nonumber\\
 &=J^2\,\tau_b(\varphi)-h(d\varphi(T),J\,\tau_b(\varphi))\,J\,d\varphi(T)\nonumber\\
 &=-\,\tau_b(\varphi)-h(d\varphi(T),J\,\tau_b(\varphi))\,J\,d\varphi(T). 
 \end{align}
 Inserting (7.15) into (7.9), 
 the left hand side of (7.9) is equal to 
 \begin{align}
 \frac{3c}{4}&\sum_{j=1}^{2n}h(d\varphi(X_j),J(\,(J\,\tau_b(\varphi))^{\top})\,)\,d\varphi(X_j)\nonumber\\
 &\qquad +\sum_{i,j=1}^{2n}h(\tau_b(\varphi),B_{\varphi}(X_i,X_j)\,B_{\varphi}(X_i,X_j)\nonumber\\
 &\qquad-\frac{2nc}{4}\,\tau_b(\varphi)
 +\frac{3c}{4}
 \,\bigg\{
 -\tau_b(\varphi)-h(d\varphi(T),J\,\tau_b(\varphi))\,Jd\varphi(T)
 \bigg\}\nonumber\\
 &=
 \sum_{i,j=1}^{2n}h(\tau_b(\varphi),B_{\varphi}(X_i,X_j))\,B_{\varphi}(X_i,X_j)\nonumber\\
 &\qquad-\frac{c(2n+3)}{4}\,\tau_b(\varphi)
 -\frac{3c}{4}\,h(d\varphi(T),J\,\tau_b(\varphi))\,J\,d\varphi(T), 
 \end{align}
 where we used (7.13) for vanishing the first term of the left hand side of (7.16). 
 \par
 Due to (7.9) and (7.16), we obtain the equivalence relation that 
 $\varphi$ is biharmonic if and only if both the equations
 \begin{align}
 (1)\qquad\qquad
 h(d\varphi(T),J\,\tau_b(\varphi))\,(J\,d\varphi(T))^{\top}=0,
 \qquad\qquad \qquad\qquad\qquad\quad
 \end{align}
 and 
 \begin{align}
 (2)\qquad 
\qquad &\sum_{i,j=1}^{2n}h(\tau_b(\varphi),B_{\varphi}(X_i,X_j))\,B_{\varphi}(X_i,X_j)
 -\frac{c(2n+3)}{4}\,\tau_b(\varphi)
 \qquad\nonumber\\
&\qquad -\frac{3c}{4}\,h(d\varphi(T),J\,\tau_b(\varphi))\,
 (J\,d\varphi(T))^{\perp}=0, 
 \end{align}
 hold. 
 \medskip\par
 {\bf 7.3} \quad For the first equation $(1)$ (7.17) is equivalent to that 
 \begin{align}
 h(d\varphi(T),J\,\tau_b(\varphi))=0\quad\mbox{or}\quad 
 (J\,d\varphi(T))^{\top}=0.
 \end{align}
 But, by (7.10), we have
 \begin{align}
 h(d\varphi(T),J\,\tau_b(\varphi))&=\big(\sum_{k=1}^{2n}H_{kk}\big)\,h(d\varphi(T),J\,\xi)\nonumber\\
 &=-\big(\sum_{k=1}^{2n}H_{kk}\big)\,
 h(J\,d\varphi(T),\xi).
 \end{align}
 By our assumption that the pseudo mean curvature 
 $\sum_{k=1}^{2n}H_{kk}\not=0$, to hold that 
 $h(d\varphi(T),J\,\tau_b(\varphi))=0$ is equivalent to 
 that 
 \begin{align}
 h(J\,d\varphi(T),\xi)=0.
 \end{align}
 And to hold that $(J\,d\varphi(T))^{\top}=0$ is equivalent to that 
 \begin{align}
 J\,d\varphi(T)=h(J\,d\varphi(T),\xi)\,\xi. 
 \end{align}
 Thus, $(1)$ (7.17) holds if and only if 
 \begin{align*}
 &\mbox{(7.21)\qquad  $h(J\,d\varphi(T),\xi)=0$, \qquad\qquad or} \\
& \mbox{(7.22)\qquad  $J\,d\varphi(T)=h(J\,d\varphi(T),\xi)\,\xi$. 
 } 
 \end{align*}
 \par
 In the case (7.21) holds, we have 
 \begin{align}
 h(d\varphi(T),J\,\tau_b(\varphi))\,(J\,d\varphi(T))^{\top}=0,
 \end{align}
which implies that (2) (7.18) turns out that 
\begin{align}
\sum_{i,j=1}^{2n}
h(\tau_b(\varphi),B_{\varphi}(X_i,X_j)\,B_{\varphi}(X_i,X_j)-\frac{c(2n+3)}{4}\,\tau_b(\varphi)=0.
\end{align}
  \par
  In the case that (7.22) holds, we have that 
  \begin{align}
  h(d\varphi(T),&J\,\tau_b(\varphi))\,(J\,d\varphi(T))^{\top}
  \nonumber\\
  &=h(d\varphi(T),\,J\,\tau_b(\varphi))\,h(J\,d\varphi(T),\,\xi)\,\xi
  \nonumber\\
  &=\big(\sum_{k=1}^{2n}H_{kk}\big)\,h(d\varphi(T),\,J\,\xi)\,h(J\,d\varphi(T),\,\xi)\,\xi
  \qquad(\mbox{by (7.10)}
  \nonumber\\
  &=-\big(\sum_{k=1}^{2n}H_{kk}\big)\,
  h(J\,d\varphi(T),\,\xi)^2\,\xi.
  \end{align}
  In the case that (7.21) holds, $(2)$ (7.18) turns out that 
  \begin{align}
  \sum_{i,j=1}^{2n}&h(\tau_b(\varphi),B_{\varphi}(X_i,X_j))-\frac{c(2n+3)}{4}\,\tau_b(\varphi)\nonumber\\
  &+\frac{3c}{4}\,\big(\sum_{k=1}^{2n}H_{kk}\big)\,h(J\,d\varphi(T),\,\xi)^2\,\xi=0. 
  \end{align}
  By inserting (7.10) into (7.24), 
  the left hand side of (7.24) is equal to 
  \begin{align}
  \sum_{i,j=1}^{2n}&h\big(\sum_{k=1}^{2n}H_{kk}\,\xi,\,H_{ij}\,\xi\big)\,H_{ij}\,\xi-\frac{c(2n+3)}{4}\,\sum_{k=1}^{2n}H_{kk}\,\xi
  \nonumber\\
  &=\big(\sum_{k=1}^{2n}H_{kk}\big)\,\bigg\{
  \sum_{i,j=1}^{2n}H_{ij}{}^2-\frac{c(2n+3)}{4}
  \bigg\}\,\xi.
  \end{align}
  $(2)$ (7.18) is equivalent to that 
  \begin{align}
   \sum_{i,j=1}^{2n}H_{ij}{}^2=\frac{c(2n+3)}{4}
  \end{align}
  by our assumption that 
  $\sum_{k=1}^{2n}H_{kk}\not=0$. 
  \par
  In the case that (7.22) holds, by inserting (7.10) into 
  (7.26), the left hand side of (7.26) is equal to 
  \begin{align}
  &\sum_{i,j=1}^{2n}h\big(\sum_{k=1}^{2n}H_{kk}\,\xi,\,H_{ij}\,\xi\big)\,H_{ij}\,\xi-\frac{c(2n+3)}{4}\,\sum_{k=1}^{2n}H_{kk}\,\xi
  \nonumber\\
  &\qquad\qquad\qquad+\frac{3c}{4}\,h(J\,d\varphi(X),\,\xi)^2\,\big(\sum_{k=1}^{2n}H_{kk}\big)\,\xi\nonumber\\
  &=\big(\sum_{k=1}^{2n}H_{kk}\big)\,\bigg\{
  \sum_{i,j=1}^{2n}H_{ij}{}^2-\frac{c(2n+3)}{4}
  +\frac{3c}{4}\,h(J\,d\varphi(T),\,\xi)^2
  \bigg\}\,\xi.
  \end{align}
Since (7.22) $J\,d\varphi(T)=h(J\,d\varphi(T),\,\xi)\,\xi$, 
we have 
\begin{align*}
h(J\,d\varphi(T),\xi)^2&=h(J\,d\varphi(T),d\,\varphi(T))\\
&=h(d\varphi(T),d\varphi(T))=g_{\theta}(T,T)=1
\end{align*} 
which implies again by our assumption 
$\sum_{k=1}^{2n}H_{kk}\not=0$, that 
$(2)$ (7.18) is equivalent to that 
\begin{align}
\sum_{i,j=1}^{2n}H_{ij}{}^2-\frac{c(2n+3)}{4}+\frac{3c}{4}=\sum_{i,j=1}^{2n}H_{ij}{}^2-\frac{n}{2}\,c=0. 
\end{align}
\par
Therefore, we obtain 
\medskip\par
\begin{th}
\quad Assume that $\varphi:\,(M,g_{\theta})\rightarrow {\mathbb P}^{n+1}(c)=(N,h)$ $(c>0)$ is an admissible isometric immersion whose pseudo mean curvature vector filed along $\varphi$ is parallel with 
non-zero pseudo mean curvature. Then, 
$\varphi$ is biharmonic if and only if 
one of the following two cases occurs: 
\par
$(1)$ \quad $h(J\,d\varphi(T),\,\xi)=0$ and 
\begin{align}
\big\Vert\, B_{\varphi}\big\vert_{H(M)\times H(M)}\,\big\Vert^2= \frac{c(2n+3)}{4},
\end{align} 
\par
$(2)$ \quad $J\,d\varphi(T)=h(J\,d\varphi(T),\xi)\,\xi$ and 
\begin{align}
\big\Vert \,B_{\varphi}\big\vert_{H(M)\times H(M)}\,\big\Vert^2= \frac{n}{2}\,c.
\end{align} 
\end{th}
\medskip\par
\section{Examples of pseudo harmonic maps and pseudo biharmonic maps}
In this section, we give some examples of pseudo biharmonic maps. 
\bigskip\par
{\em Example 8.1.} \quad Let 
$
(M^{2n+1},g_{\theta})=S^{2n+1}(r)$ be the sphere of radius 
$r$ $(0<r<1)$ which is 
embedded in the unit sphere $S^{2n+2}(1)$,\,i.e., 
\par\noindent 
the natural embedding $\varphi:\,S^{2n+1}(r)\rightarrow S^{2n+2}(1)$ is given by
$$ 
\varphi:\,
S^{2n+1}(r)\ni
x'=(x_1,x_2,\cdots,x_{2n+2})
\mapsto 
(x',\sqrt{1-r^2})\in S^{2n+2}(1).
$$
This $\varphi$ is a standard isometric with constant principal curvature $\lambda_1=\cot\,[\,\cos^{-1}t\,]$, $(-1<t<1)$, with 
the multiplicity $m_1=\dim M=2n+1$.  
\par
Due to Theorem 6.2, 
it is {\em pseudo biharmonic} if and only if 
\begin{align}
(\lambda_1)^2\,\times\,2n=2n\qquad
&\Leftrightarrow\quad 
\lambda_1=\cot\,[\,\cos^{-1}t] 
=\pm1.
\nonumber\\
&\Leftrightarrow\quad 
\displaystyle t=\cos\left(\pm\frac{\pi}{4}\right)=\frac{1}{\sqrt{2}}.
\end{align}
This is just the example which is biharmonic but not minimal given by 
C. Oniciuc (\cite{O}). Note that 
$\varphi:\,S^{2n+1}(r)\rightarrow S^{2n+2}(1)$ is 
{\em pseudo harmonic} if and only if 
\begin{align}
\mbox{\rm Trace}
(B_{\varphi}
\vert_{H(M)\times H(M)}
)=0
&\quad \Leftrightarrow\quad \lambda_1=0
\nonumber\\
&\quad \Leftrightarrow\quad 
t=\cos\,\left(\,\frac{\pi}{2}\,\right)=1.
\end{align}
This $t=1$ gives a great hypersphere which is also minimal.
\bigskip\par
{\em Example 8.2.}\quad 
Let 
the Hopf fibration
$\pi:\,S^{2n+3}(1)\rightarrow {\mathbb P}^{n+1}(4)$, 
and,  
let $\widehat{M}:=S^{1}(\cos \,u)\times 
S^{2n+1}(\sin \,u)\,\subset\,S^{2n+3}(1)$ 
$\left( 0<u<\frac{\pi}{2} \right)$. 
Then, we have 
$\varphi:\,M^{2n+1}
=\pi(\widehat{M})\,\subset {\mathbb P}^{n+1}(4)$ 
which is a homogeneous real hypersurface of ${\mathbb P}^{n+1}(4)$ of type $A_1$ in the table of R. Takagi (\cite{RT})
whose principal curvatures and their multiplicities are given as follows (\cite{RT}): 
\begin{align}
\left\{
\begin{aligned}
&\lambda_1=\cot\,u,\qquad
\mbox{\rm multiplicity}\,\,m_1=2n,\\
&\lambda_2=2\,\cot(\,2u\,),\qquad\mbox{\rm multiplicity}\,\,m_2=1.
\end{aligned}
\right.
\end{align}
Since $2\,\cot(\,2u\,)=\cot\,u-\tan\,u$, 
the mean curvature $H$ and $\Vert B_{\varphi}\Vert^2$ are given by 
\begin{align}
H&=\frac{1}{2n+1}\,\left\{
(2n+1)\cot\,u-\tan\,u
\right\},\\
\Vert\,B_{\varphi}\,\Vert^2&=m_1\,\lambda_1{}^2+m_2\,\lambda_2{}^2=\tan^2u+(2n+1)\,\cot^2u-2.
\end{align}
R. Takagi showed (\cite{RT}) to this example, 
that 
$\varphi:\,M^{2n+1}\rightarrow {\mathbb P}^{n+1}(4)$ 
is the geodesic sphere $S^{2n+1}$, 
and $J(-\xi)$ is the mean curvature vector of the principal curvature 
$\lambda_2$ (cf. Remark 1.1 in \cite{RT}, p. 48), 
where $\xi$ is a unit normal vector field along $\varphi$. 
\medskip\par
In the case $(1)$ of Theorem 7.1, i.e., 
$(M^{2n+1},g_{\theta})=(S^{2n+1},g_{\theta})$ 
is a strictly pseudoconvex $CR$ manifold and 
$J\,d\varphi(T)$ is tangential, 
we have 
$$
0=h(J\,d\varphi(T),\xi)=h(J^2\,d\varphi(T),J\,\xi)=h(d\varphi(T),J(-\xi)), 
$$
and 
$h(d\varphi(T),d\varphi(H(M)))=0$. 
Then, the principal curvature vector field $J(-\xi)$ 
with principal curvature $\lambda_2=2\,\cot(2u)$ 
coincides with $d\varphi(X)$ for some $X\in H(M)$. 
Since
$$
\Vert X\Vert=\Vert d\varphi(X)\Vert=\Vert J(-\xi)\Vert=\Vert \xi\Vert=1,
$$
we can choose an orthonormal basis $\{X_i\}_{i=1}^{2n}$ of $H(M)$ in such a way $X_1=X$. 
Then, 
$\{d\varphi(T),d\varphi(X_2),\cdots,d\varphi(X_{2n})\}$ give principal curvature vector fields along $\varphi$ 
with principal curvature $\lambda_1=\cot\,u$. 
Then, 
\begin{align}
\tau_b(\varphi)&=\sum_{i=1}^{2n}B_{\varphi}(X_i,X_i)
=2\,\cot (2u)+(2n-1)\,\cot\,u\nonumber\\
&=2n\,\cot u-\tan u. 
\end{align}
Therefore, $\varphi$ is {\em pseudo harmonic} if and only if 
\begin{align}
\tau_b(\varphi)=0\qquad\Leftrightarrow\quad 
\tan u=\sqrt{2n}. 
\end{align}
\par
By Theorem 7.1, $(1)$, $\varphi$ is pseudo biharmonic if and only if 
\begin{align}
\Vert B_{\varphi}\vert_{H(M)\times H(M)}\Vert^2=\frac{c(2n+3)}{4}=2n+3.
\end{align}
Since the left hand side of (8.8) coincides with 
\begin{align}
\Vert B_{\varphi}\vert_{H(M)\times H(M)}\Vert^2&=(2\,\cot (2u))^2+(2n-1)\,\cot^2u\nonumber\\
&=(\cot u-\tan u)^2+(2n-1)\cot^2u\nonumber\\
&=2n\cot^2u-2+\tan^2u, 
\end{align}
we have that $(8.8)$ holds if and only if  
\begin{align}
2n\,\cot^2u+\tan^2u=2n+5
\qquad\Leftrightarrow\quad  x^2-(2n+5)x+2n=0,
\end{align}
where $x=\tan^2u$. 
Therefore, $\varphi$ is {\em pseudo biharmonic} 
if and only if $\tan u$ is $\sqrt{\alpha}$ or $\sqrt{\beta}$, 
where $\alpha$ and $\beta$ are positive roots of 
(8.10). 
\medskip\par
In the case $(2)$ of Theorem 7.1, i.e.,  
$(M^{2n+1},g_{\theta})=S^{2n+1},g_{\theta})$ is a strictly pseudoconvex 
$CR$ manifold, and $J\,d\varphi(T)$ is normal, i.e., 
$J\,d\varphi(T)=h(J\,d\varphi(T)\,\xi)\,\xi$. Then, we have that 
$$0\not=d\varphi(T)=h(d\varphi(T),J(-\xi))\,J(-\xi).$$
And $J(-\xi)$ is the principal curvature vector field along $\varphi$ 
with the principal curvature $\lambda_2$, and 
$d\varphi(H(M))$ is the space spanned by 
the principal curvature vectors along $\varphi$ with the principal curvature 
$\lambda_1$ since $h(d\varphi(H(M)), J(-\xi))=0$.  Then 
the pseudo tension field $\tau_b(\varphi)$ is given by 
\begin{align}
\tau_b(\varphi)
=\sum_{i=1}^{2n}B_{\varphi}(X_i,X_i)=(2n\,\cot\,u)\,\xi\not=0, 
\end{align}
so that $\varphi$ is {\em not pseudo harmonic}.  
Due to the case (2) of Theorem 7.1 that $J\,d\varphi(T)$ is normal, 
$\varphi$ is pseudo biharmonic if and only if 
\begin{align}
2n&=\Vert B_{\varphi}\vert_{H(M)\times H(M)}\Vert^2\nonumber\\
&=m_1\lambda_1{}^2\nonumber\\
&=2n\cot^2u
\end{align}
occurs. 
Thus, we obtain  
\begin{align}
\tan^2u=1. 
\end{align}
Therefore, if 
$\tan\,u=1$ ($u=\frac{\pi}{4}$), then 
the corresponding isometric immersion $\varphi:\,(M^{2n+1},g_{\theta})\rightarrow {\mathbb P}^{n+1}(4)$ is {\em pseudo biharmonic}, 
but {\em not pseudo harmonic}.   
\medskip\par
\begin{rem}
Let us recall our previous work $($\cite{IIU2}, \cite{IIU}$)$ 
that $\varphi:\,(M^{2n+1},g_{\theta})\rightarrow {\mathbb P}^{n+1}(4)$ is {\em biharmonic} if and only if 
\begin{align}
\Vert B_{\varphi}\Vert^2&=\tan^2u+(2n+1)\,\cot^2u-2=\frac{n+2}{2}\,4\nonumber\\
&\Leftrightarrow\,\,
x^2-2(n+3)x+2n+1=0,\qquad (x=\tan^2u).
\end{align}
The equation (8.14) has two positive solutions $\alpha,\,\beta$, 
and 
if we put $\tan \,u=\sqrt{\alpha}$ or $\sqrt{\beta}$ $(0<u<\frac{\pi}{2})$, 
 then $\varphi:\,(M^{2n+1},g_{\theta})\rightarrow {\mathbb P}^{n+1}(4)$ is {\em biharmonic}, and vice versa. 
 Since the mean curvature is given by (8.4), 
 $\varphi:\,(M^{2n+1},g_{\theta})\rightarrow {\mathbb P}^{n+1}(4)$ is {\em harmonic} (i.e., {\em minimal}) 
 if and only if 
 $\tan\,u=\sqrt{2n+1}$ $(0<u<\frac{\pi}{2})$. 
 \end{rem}
\medskip\par

\end{document}